\documentclass[11pt]{amsart}
\usepackage{amsmath}

\def\ol{\overline}
\def\e{\epsilon}
\def\lf{\left}
\def\ri{\right}
\def\a{{\alpha}}
\def\b{{\beta}}

\def\bb{{\bar\b}}

\def\abb{{\a\bb}}

\def\wt{\widetilde}

\def\p{\partial}
\newcommand\rr{{\mathbb R}}
\newcommand\ce{{\mathbb C}}
\newcommand\C{{\mathbb C}}

\def\ii{\sqrt{-1}}
\def\jbar{{\bar\jmath}}

\def\ttg{\tilde{g}}

\def\K{K\"ahler }

\def\KRF{K\"ahler-Ricci flow }
\def\KRS{K\"ahler-Ricci soliton }
\def\KRS{K\"ahler-Ricci soliton }

\def\be{\begin{equation}}
\def\ee{\end{equation}}
\def\ol{\overline}
\def\lf{\left}
\def\ri{\right}
\def\a{{\alpha}}
\def\b{{\beta}}

\def\e{\epsilon}

\def\ijb{{i\jbar}}

\def\bb{{\bar\b}}

\def\Ric{\text{Ric}}

\def\abb{{\a\bb}}

\def\wt{\widetilde}

\def\p{\partial}

\def\C{\Bbb C}
\def\cn{\Bbb C^n}
\def\wh{\widehat}

\def\wt{\widetilde}

\def\p{\partial}

\def\C{\Bbb C}
\def\ii{\sqrt{-1}}
\def\KRF{K\"ahler-Ricci flow }
\def\ttg{\wt{g}}
\def\ttR{\wt  R}
\newtheorem{thm}{Theorem}[section]

\newtheorem{lem}{Lemma}[section]
\newtheorem{prop}{Proposition}[section]
\newtheorem{cor}{Corollary}[section]
\theoremstyle{definition}

\theoremstyle{remark}
\newtheorem{rem}{Remark}
\numberwithin{equation}{section}
 
\begin{document}

\title{On the complex structure of K\"ahler manifolds with nonnegative
 curvature}

%    Information for first author
\author{Albert Chau}
%    Address of record for the research reported here
\address{Harvard University, Department of Mathematics,
  One Oxford Street, Cambridge, MA 02138, USA}
\email{chau@math.harvard.edu}
%    Current address
%\curraddr{}
%    \thanks will become a 1st page footnote.
%\thanks{}

%    Information for second author
\author{Luen-Fai Tam$^1$}

\thanks{$^1$Research
partially supported by Earmarked Grant of Hong Kong \#CUHK4032/02P}

\address{Department of Mathematics, The Chinese University of Hong Kong,
Shatin, Hong Kong, China.} \email{lftam@math.cuhk.edu.hk}

%    General info
\renewcommand{\subjclassname}{%
  \textup{2000} Mathematics Subject Classification}
% Global differential geometry
% 53C21 Methods of Riemannian geometry, including PDE methods;
%   curvature restrictions [See also 58J60]
% 53C44 Geometric evolution equations (mean curvature flow)
% 53C55 Hermitian and K\"ahlerian manifolds [See also 32Cxx]
% Qualitative properties of solutions
% 35B35 Stability, boundedness
% Partial differential equations on manifolds; differential operators
% 58J37 Perturbations; asymptotics
% Systems theory; control - Stability
% 93D05 Lyapunov and other classical stabilities
%       (Lagrange, Poisson, $L^p, l^p$, etc.)

\date{April 2005.}

%\dedicatory{}

%\keywords{}

\begin{abstract} We study the asymptotic behavior of the K\"ahler-Ricci
  flow on \K manifolds of nonnegative holomorphic bisectional curvature. Using these results we prove
   that a complete noncompact
  K\"ahler manifold with nonnegative bounded holomorphic bisectional curvature
  and maximal volume growth is biholomorphic to complex Euclidean
  space $\C^n$. We also show that the volume growth condition can be
removed if we assume $(M, g)$ has average quadratic scalar
curvature decay (see Theorem \ref{shitamthm}) and positive curvature operator.
\end{abstract}

\maketitle

\markboth{Albert Chau and Luen-Fai Tam} {Complex structure of K\"ahler manifolds with nonnegative
 curvature}

\section{Introduction}

The classical uniformization theorem says that a simply connected Riemann surface
is either the Riemann sphere, the open unit disk or the complex plane.  On the
other hand, there is a close relation between the complex structure and the
geometry of a Riemann surface. An important case of this is that a complete
noncompact Riemannian surface with positive Gaussian curvature is necessarily
conformally equivalent to the complex plane. In higher dimensions, there is a
long standing conjecture predicting similar results.  In its most general form,
the conjecture is due to Yau \cite{Y}, and it states: {\it A complete noncompact
K\"ahler manifold with positive holomorphic bisectional curvature is
biholomorphic to $\ce^{n}$}. In fact, the conjecture is part of a program
proposed by Yau \cite{Y} in 1974 to study complex manifolds of parabolic type.

The first result supporting this conjecture was by Mok-Siu-Yau
\cite{MSY}. There, the authors proved that if $M^n$ is a complete noncompact
K\"ahler manifold with nonnegative bisectional curvature, maximal
volume growth and faster than
quadratic scalar curvature decay, then $M$ is isometrically biholomorphic to $\cn$.  Later,
Mok \cite{M1} proved that if $M$ has positive bisectional curvature,
maximal volume growth and quadratic scalar curvature decay,
then $M$ is an affine algebraic variety.  As a consequence, if $n=2$
and the sectional curvature is positive, then $M$ is biholomorphic to
$\Bbb C^2$ by a result of  Ramanujan \cite{R}. In this case, dimension
2, it is known that the condition
on the sectional curvature can be relaxed and the decay of the scalar curvature can also be removed, see \cite{CTZ, CZ3,
  Ni2}. In higher dimensions and in general, the conjecture is still
very open, and until now, this has been so even if $M$ is
assumed to have bounded curvature and maximal
volume growth.  In this paper (Corollary \ref{s1c1}) we show that the
conjecture is true in
all dimensions provided $M$ has bounded curvature and maximal volume growth.

In his thesis \cite{Sh}, Shi used the following Ricci flow of Hamilton \cite{Ha}
 to better understand the uniformization conjecture
in the case of $(M, g)$ as in Mok's paper \cite{M1}:

\begin{equation}\label{s1e1}
\begin{split}
\frac{\p }{\p t}\tilde{g}_{i\jbar}(x,t)&=-\tilde{R}_{i\jbar}(x,t)\\
\tilde{g}_{i\jbar}(x,0)&=\tilde{g}_{i\jbar}(x).\\
\end{split}
\end{equation}

On a \K manifold, (\ref{s1e1}) is referred to as the \KRF.  In \cite{Sh,Sh0}, Shi
obtained several important results for this flow including short time existence
for general solutions, and long time existence together with many useful
estimates in the above case; see Theorem \ref{shitamthm} for more details.
Although the results in \cite{Sh} did not actually prove uniformization in this
case\footnote{This was observed later on in \cite{Chen}.  Also see \cite{CT}.},
their importance has remained fundamental to the field; in particular, in the
above mentioned works \cite{CTZ}, \cite{CZ3}, \cite{Ni2} as well as the present
paper.

 In this paper, by studying the asymptotic behavior of the K\"ahler-Ricci flow (\ref{s1e1})
 in more detail, we will prove the following uniformization theorem:

\begin{thm}\label{s1t1}
Let $(M^n,\wt g)$ be a
complete noncompact K\"ahler manifold with nonnegative and bounded holomorphic
bisectional curvature.  Suppose
\begin{enumerate}
\item[(i)] $Vol(B(p, r)) \geq C_1r^{2n};\hspace{11pt} \forall r\in [0,
  \infty)$  for some   $p \in M$,
\item[(ii)]$\frac{1}{V_x(r)}\int_{B_x(r)}R\le \frac{C_2}{1+r^2}$
for all $x\in M$ and for all $r>0$,
\end{enumerate}
for some positive constants $C_1, C_2$.  Then $M$ is biholomorphic to
$\C^n$.  Moreover, condition (i) can be removed if $M$ has positive curvature operator.
\end{thm}
The assumption is a little bit weaker than in \cite{Sh}. Here we only assume the
holomorphic bisectional curvature is nonnegative.  

In \cite{Y}, Yau conjectured that (i) actually implies (ii). This has been confirmed by  Chen-Tang-Zhu \cite{CTZ}  for the case of dimension 2, and Chen-Zhu \cite{CZ3} for higher dimensions under the additional  condition that the curvature operator is nonnegative and recently by Ni \cite{Ni2} for all dimensions. Hence we have:  
\begin{cor}\label{s1c1}
Let $(M^n,\wt g)$ be a complete noncompact K\"ahler manifold
with nonnegative and bounded holomorphic bisectional
curvature and maximal volume growth.  Then $M$ is biholomorphic to $\C^n$.
\end{cor}

Also, if  we only assume that the curvature operator is nonnegative and if the scalar curvature is still assumed to satisfy (ii), then one can prove that   the universal cover of $M$ is
biholomorphic to $\C^n.$

In order to prove  Theorem \ref{s1t1}, we will first obtain some   results
on the long time behavior of the \KRF (\ref{s1e1}).  We believe these
results may be of independent interest.  It will be more convenient to consider the   normalized K\"ahler-Ricci flow
 \begin{equation}\label{s1e2}
\frac{\p }{\p t}g(t)=-Rc(t)-g(t)
\end{equation}
where $g(t)=e^{-t}\wt g(e^t)$ and $Rc(t)$ is the Ricci curvature of $g(t)$. Under the assumptions of Theorem \ref{s1t1}, we have:
\begin{thm}\label{s1t2} Let $(M^n,\ttg)$ be as in Theorem
  \ref{s1t1} with maximal volume growth or with positive curvature operator
  and let $g(x,t)$ be as in (\ref{s1e2}). Let $p\in M$ be any point.
Then the eigenvalues of $Rc(p,t)$ with respect to $g(p,t)$ will converge as $t\to\infty$. Moreover, if $\mu_1>\mu_2>\dots>\mu_l$ are the distinct limits of the eigenvalues, then
   $V=T_p^{(1,0)}(M)$ can be decomposed orthogonally with respect to $g(0)$ as
 $V_1\oplus \cdot\cdot\cdot\oplus V_l$ so that the following are true:
\begin{enumerate}
\item[(i)] If $v$  is a nonzero vector in  $V_i$ for some $1\le i\le l$, then
 $$\lim_{t \to \infty} Rc(v(t),\bar v(t))=\mu_i$$ and thus
$$\lim_{t \to \infty} \frac{1}{t}\log \frac{|v|_{g(t)}^2}{|v|_{g(0)}^2}=-\mu_i-1.$$  Moreover, both convergences are uniform over all
 $v\in V_i\setminus\{0\}$. Here $v(t)=v/|v|_{g(t)}$.

\item[(ii)] For $1\le i, j\le l$ and for nonzero vectors  $v \in V_i$ and $w \in V_j$ where $i\neq j$, $\lim_{t\to
\infty}\langle v(t),w(t)\rangle_t=0$ and the convergence is uniform over all such nonzero vectors $v, w$.
\item[(iii)] $\dim_\C(V_i)=n_{i}-n_{i-1}$  for each $i$.
\item[(iv)] $$\sum_{i=1}^l(-\mu_i-1)\dim_\C
V_i=\lim_{t\to\infty} \frac1t\log \frac{\det(g_{i\bar j}(t))}{\det ({g}_{i\bar
j}(0)}.$$ \end{enumerate}
\end{thm}
The theorem says that at $p$, $(M,g(t))$ behaves like a gradient
K\"ahler-Ricci soliton of expanding type, see Proposition \ref{s3p2}
for more details. Moreover, conclusions (i) and (ii) mean in some
sense that $Rc(p,t)$ can be `simultaneously digonalized' near
infinity. Conditions (i), (iii) and (iv) basically say  that $g(t)$ is
Lyapunov regular in the sense of dynamical systems, see \cite{BP}.

On the other hand, by Theorem \ref{shitamthm}, we can construct biholomorphisms
from a sequence of open sets which exhaust $M$ onto a fixed ball in $\cn$.  By
identifying these open sets, the results in Theorem \ref{s1t2} can be interpreted
in terms of the dynamics of a randomly iterated sequence of biholomorphisms as in
\cite{JV}. Using the results of Theorem \ref{s1t2} in this setting, and using
techniques developed by Rosay-Rudin \cite{RR} and Jonsson-Varolin \cite{JV}, we
then proceed to assemble these biholomorphisms into a  global biholomorphism from
$M$ to $\cn$.

The paper is organized as follows.  In \S2 we review the main results
on the \KRF (\ref{krf}) which we use later.  In \S3 and \S4 we study
the asymptotic behavior of the \KRF on $M$ as $t \to \infty$.  The
focus of \S3 will primarily be on the global asymptotics of the \KRF on $M$ while that of
\S4 will be purely local.  We believe that these asymptotics should be of
independent interest to the study of the \KRF.  Finally, in \S5 we
will prove Theorem \ref{s1t1} and its corollaries.

The authors would like to thank Richard Hamilton and S.T.Yau for helpful
discussions

\section{the K\"ahler Ricci flow}

In this section, we will collect some known results on \KRF which will be used in this work. Recall that on a complete noncompact \K manifold
$(M^n,\tilde{g}_\ijb(x))$,  the \KRF equation is:
\begin{equation}\label{krf}
\begin{split}
\frac{\p }{\p t}\tilde{g}_{i\jbar}(x,t)&=-\tilde{R}_{i\jbar}(x,t)\\
\tilde{g}_{i\jbar}(x,0)&=\tilde{g}_{i\jbar}(x).\\
\end{split}
\end{equation}

\begin{thm}\label{shitamthm} Let $(M^n,\tilde{g})$ be a complete noncompact
  K\"ahler manifold with bounded nonnegative holomorphic bisectional
  curvature.  Suppose there is a constant $C>0$ such that its scalar curvature satisfies
\begin{equation}\label{main1e1}
\frac1{V_x(r)}\int_{B_x(r)}\tilde{R} \hspace{3pt}dV_g\le \frac{C}{1+r^2}
\end{equation}
 for all $x\in M$ and for all $r>0$ where $\tilde{R}$ is the scalar
curvature. Then the \KRF (\ref{krf}) has a long time solution
$\tilde{g}_\abb(x,t)$ on $M\times [0, \infty)$.  Moreover, the following are
true:
\begin{enumerate}
\item[(i)] For any $t\ge 0$, $\tilde{g}(x,t)$ is K\"ahler with nonnegative
holomorphic bisectional curvature.
 \item [(ii)] For any integer $m\ge 0$,
there is a constant $C_1$ depending only on $m$ and the initial metric such that
$$
||\nabla^m \tilde{R}m||^2(x,t)\le \frac{C_1}{t^{2+m}},
$$
for all $x\in M$ and for all $t\ge 0$, where $\nabla $ is the covariant
derivative with respect to $\tilde{g}(t)$ and the norm is also taken in
$\tilde{g}(t)$.
\item [(iii)] If in addition $(M,\wt g(0))$ has either maximum volume growth
or
  positive curvature operator, then there exists a constant $C_2>0$ depending
  only on  the initial metric such that  the injectivity radius of $\tilde{g}(t)$ is bounded below by $C_2t^{1/2}$ for all $t\ge1$.
\end{enumerate}
\end{thm}
 \begin{proof} (i) and (ii) are mainly obtained by   Shi \cite{Sh0, Sh, Sh2} (also see \cite{NT}). To prove (iii), suppose $\wt g(0)$ has positive curvature operator. 
 Then by \cite{Ha4} we know that positive curvature operator is preserved
 under (\ref{krf}), and thus $g(t)$ has positive sectional curvature at
every time $t$. From this and the estimates in (ii), we can conclude by the
results in \cite{GM} that (iii) is true in the case of positive curvature operator
(See \cite{CKL} p. 14 for a description of how to prove this). In the
 case of maximal volume growth,
(iii) has been observed in \cite{Chen}. In fact,
if $Vol_0(B_x(r))\ge Cr^{2n}$ for some $C>0$ for the initial metric, then we also
have $Vol_t(B_x(r))\ge Cr^{2n}$ for the metric $g(t)$ for all $t\ge0$
with the same constant $C$, see \cite{Chen} for example. Combining
this with the curvature estimates (ii) and the injectivity radius
estimates in \cite{CGT}, (iii) follows in this case.
\end{proof}
We now consider the following normalization of (\ref{krf}):
\begin{equation}\label{krfn}
\frac{\p }{\p t}g_{i\jbar}(x,t)=-R_{i\jbar}(x,t)-g_{i\jbar}(x,t).
\end{equation}
It is easy to verify that if $\tilde{g}(x, t)$ solves (\ref{krf}), then
\begin{equation}\label{krfns}
g(x, t)=e^{-t}\tilde{g}(x, e^{t})
\end{equation}
is a solution to (\ref{krfn}).  Thus for $\tilde{g}(x, t)$ as in Theorem 2.1,
$g(x, t)$ in (\ref{krfns}) is defined for $-\infty<t<\infty$. Note that
$\lim_{t\to-\infty}g(x,t)=\wt g(x)$ which is the initial data of
(\ref{krf}). The results in Theorem \ref{shitamthm} can be translated
to the following results for a solution to (\ref{krfn}):

\begin{cor}\label{pinchcor}
Let $\tilde{g}(x, t)$ be as in Theorem \ref{shitamthm} and let $g(x, t)$ be given
by (\ref{krfns}). Then the following are true:
\begin{enumerate}
\item[(i)] For any $-\infty<t<\infty$, $g(x,t)$ is K\"ahler with nonnegative
holomorphic bisectional curvature.
 \item [(ii)] For any integer $m\ge 0$,
there is a constant $C_1$ depending only on $m$ and the initial metric such that
$$
||\nabla^m  {R}m||^2(x,t)\le {C_1},
$$
for all $x\in M$ and for all $t\ge 0$, where $\nabla $ is the covariant
derivative with respect to ${g}(t)$ and the norm is also taken in
$ {g}(t)$.
\item [(iii)] If in addition $(M, \wt g(0))$ has either maximum volume growth
or positive curvature operator, then there exists a constant $C_2>0$ depending
  only on  the initial metric such that  the injectivity radius of $ {g}(t)$ is bounded below by $C_2$ for all $t\ge0$.
\end{enumerate}
\end{cor}

We shall need the follwing.

\begin{prop}\label{holomorphicoord} Let $(M^n,g)$ be a complete K\"ahler manifold with nonnegative holomorphic bisectional curvature such that $|Rm|+|\nabla Rm| \le C_1$ and the injectivity of $M$ is larger than $r_0$. Then there exist positive constants $r_1,r_2$ and $C_2$ depending only on $C_1$, $r_0$ and $n$ such that for each $p\in M$, there is a  holomorphic map $\Phi$ from the Euclidean ball $\wh B_0(r_1)$ at the origin of $\cn$ to $M$ satisfying the following:
\begin{enumerate}
\item [(i)]$\Phi$ is a biholomorphism from $\wh B_0(r_1)$ onto its image;
\item[(ii)] $\Phi(0)=p$;
\item[(iii)] $\Phi^*(g)(0)=g_\e$;
\item[(iv)] $\frac{1}{r_2}g_{\epsilon}\leq\Phi^*(g) \leq  r_2g_{\epsilon}$ in $\wh B(0, r_1)$.
\end{enumerate}
where $g_\e$ is the standard metric on $\cn$.
\end{prop}
\begin{proof} This is in fact a special case of  Proposition 1.2 in \cite{TY}, see also \cite{Sh2, Chen}. For the sake of completeness, we sketch the proof as follows.

 By the assumption on the injectivity radius, let $x_1,\dots,x_{2n}$ be normal coordinates on $B_{p}(r_0)$ so that if $z_i=x_i+\sqrt{-1}x_{n+i}$ are standard complex coordinates of $\cn$, then $\frac{\p}{\p z_i}$ form a basis for $T_{p}(M)$ at $p$.  Hence there is a diffeomorphism $F$ from $B_{p}(r_0)$ onto $\wh B_0(r_0)$   such that $F(p)=0$ and $dF\circ  J=\wh J\circ dF$ at $0$ where $\wh J$ is the standard complex structure on $\cn$ and $J$ is the complex structure of $M$. By \cite{Ha2},  the components of the metric $g$  with respect to  coordinates $x_i$ satisfies
$$
|\delta_{ij}-g_{ij}|\le C_2|x|^2, \ \frac12\delta_{ij}\le g_{ij}\le 2\delta_{ij},
$$
$$
|\frac{\p^2}{\p x_k\p x_l}g_{ij}|\le C_2
$$
and
$$
|\frac{\p }{\p x_k}g_{ij}|(x)\le C_2|x|
$$
in $B_{p}(r_1)$ for some positive constants $r_1$, $C_2$ depending only on $C_1$, $r_0$  and $n$. Here $|x|^2=\sum_{i}(x_i)^2$. In the following $C_i$'s and $r_i$'s always denote  positive constants depending only on $C_1$, $r_0$ and $n$.
Hence if $r_1$ small enough,   $\sqrt{-1}\p\ol\p\log \rho^2\ge -C_3\omega$ and the eigenvalues of the Hessian of $\rho^2$ are bounded below by $C_4$. Here $\rho$ is the distance from $p$ and $\omega$ is the K\"ahler form.  One can   prove that $|J-\wh J|\le C_5\rho^2$, where we also denote the pull back of $\wh J$ under $F$ with $\wh J$, see \cite{Chen}. The $i$-th component $z_i=x_i+\sqrt{-1}x_{n+i}$ of the map $F$ when considered as a map from $B_{p}(r_0)$ to $\cn$ satisfies

\begin{equation}\label{dbarbound}
|\ol\p z_i|\le C_6\rho^2.
\end{equation}
As in \cite{Ni}, by Corollary 5.3 in \cite{Demailly}, using the weight function $\varphi=(n+2)\log \rho^2+C_7\rho^2$ for some $C_7$ so that $\sqrt{-1}\p\ol \p\varphi\ge C_8\omega$, one can solve $\ol\p u_i=\ol \p z_i$ in $B_{p}(r_1)$ with
\begin{equation}\label{L2estimate}
\int_{B_{p}(r_1)}|u_i|^2e^{-\varphi}\le \frac1{C_8}\int_{B_{p}(r_1)}|\ol \p z_i|^2e^{-\varphi}\le C_9
\end{equation}
for some $C_9$. Here we have used the fact that $Ric\ge0$ and (\ref{dbarbound}). From this, it is easy to see that $u_i(p)=0$ and $du_i(p)=0$. Moreover, from the fact that $z_i-u_i$ is holomorphic one can prove that on $B_{p}(r_1/2)$,
$$
|u_i|+|\nabla u_i|+|\nabla^2 u_i|\le C_{10}
$$
by (\ref{L2estimate}), mean value inequality, gradient estimates and
Schauder estimates. Hence we have  $|\nabla u_i|\le C_{11}\rho$ and
$|u_i|\le C_{11}\rho^2$. So the map $\Phi$   given by
$\Phi^{-1}=(z_1-u_1,\dots,z_n-u_n)$ will satisfy the conditions in the
proposition if $r_1$ is small enough and  
$r_2$ is large enough.
\end{proof}
 Using this and Corollary \ref{pinchcor}, we have the following (also
 see \cite{TY, Sh2}).

 \begin{cor}\label{Assumption} Let $(M^n, \wt g(0))$ and $g(x,t)$ be as in Corollary \ref{pinchcor}
such that $(M,\wt g(0))$ has either maximum volume growth or positive
  curvature operator. Let $p\in M$ be a fixed point. Then there are constants $r_1$ and $r_2$ depending only on the initial metric such
that for every $t>0$ there exists a  holomorphic map $\Phi_t:\wh B_0( r_1)\subset \cn
\to M$ satisfying:
\begin{enumerate}
\item [(i)]$\Phi_t$ is a biholomorphism from $\wh B_0(r_1)$ onto its image;
\item[(ii)] $\Phi_t(0)=p$; 
\item[(iii)] $\Phi_t^*(g(t))(0)=g_\e$; \item[(iv)]
$\frac{1}{r_2}g_{\epsilon}\leq\Phi_t^*(g(t)) \leq  r_2g_{\epsilon}$ in $\wh B_0(r_1)$;
\end{enumerate}
where $g_\e$ is the standard metric on $\cn$, and $\wh B_0(r_1)$ is the Euclidean ball of radius $r_1$ with center at the origin in $\cn$. Moreover, the following are true:
\begin{enumerate}
\item[(v)] For any
$t_k\to\infty$ and for any $0<r<r_1$, the family $\{\Phi_{t_k}(\wh
B_0(r))\}_{k\ge 1}$ exhausts $M$ and hence  $M$ is simply
connected.
\item[(vi)] If  $T$ is large enough,
then $F_{i+1}=\Phi_{(i+1)T}^{-1}\circ \Phi_{iT}$ maps $\wh B_0(r_1)$ into
$\wh B_0(r_1)$ for each $i$, and there is $0<\delta<1$, $0<a<b<1$ such that 
 $$
|F_{i+1}(z)|\le \delta |z|
$$
for all $z\in \wh B_0(r_1)$, and
$$
a|v|\le |F'_{i+1}(0)(v)|\le b|v|
$$
for all $v$ for all $i$. 
\end{enumerate} 
\end{cor}
\begin{proof} (i)--(iv) follows immediately from Proposition \ref{holomorphicoord} and Corollary \ref{pinchcor}. To prove (v), observe that $B_p^t(r/r_2)\subset \Phi_t(\wh
B_0(r))$ by (i) and (iv), where $B_p^t(R)$ is the geodesic ball of radius $R$ with respect to $g(t)$ with center at $p$. On the other hand, by (\ref{krfn}),  $g(t)\le e^{-t/2}g(0)$ and so $B_p^0(R)\subset B_p^t(e^{-t}R)$. From this it is easy to see that (v) is true.

To prove (vi), let $v$ be a $(1,0)$ vector on $M$ and denote $|v|_t$ to be the length of $v$ with respect to $g(t)$. By (\ref{krfn}) and Corollary \ref{pinchcor} 
\begin{equation}
\begin{split}
 -|v|^2_t&\ge
\frac{d}{dt}|v|^2_t\\
&=-Rc_{\ttg}(v, v)-\ttg(v, v)\\
&\geq -  C_1\ttg(v, v)-\ttg(v, v)\\
&\geq -(  C_1+1)|v|^2_t
\end{split}
\end{equation}
for some constant  $C_1>0$ which is independent of $v$ and $t$. 
 Hence for any $T>0$ and $i\ge 1$,
\begin{equation}\label{uniformshrink11}
e^{-T}\ge \frac{|v|^2_{(i+1)T}}{|v|^2_{iT}}\ge e^{-( C_1+1)T}.
\end{equation}
Since
$$\Phi_{iT}(\wh B_0(r_1))\subset B^{iT}_{p}(r_2r_1)\subset B^{(i+1)T}_{p}(e^{-T/2}r_2r_1),$$ 
and $\Phi_{(i+1)T}(\wh B_0(r_1))\supset B^{(i+1)T}_{p}(r_1/r_2)$.  
Hence $F_{i+1}$ is defined on $\wh B_0(r_1)$ and $F_{i+1}(\wh B_0(r_1))\subset \wh B_0(r_1)$ if $T$ is large enough. From (iv) and (\ref{uniformshrink11}), it is easy to see and there is $0<\delta<1$,   such that 
 $$
|F_{i+1}(z)|\le \delta |z|
$$
for all $z\in \wh B_0(r_1)$ for all $i$ if $T$ is large. From (ii), (iii) and  (\ref{uniformshrink11}), we can also find $0<a<b<1$ such that 
$$
a|v|\le |F'_{i+1}(0)(v)|\le b|v|
$$
for all $v$ and for all $i$. This completes the proof of the corollary.
\end{proof}

 In \S5, we
will use the maps
$\Phi_t$ to construct a biholomorphism from $M$ to $\cn$.

\section{Asymptotic behavior of K\"ahler Ricci flow (I)}
Let $(M^n,\tilde{g}_\ijb(x))$ be as in Theorem \ref{shitamthm} satisfying (\ref{main1e1}).
Let
$\ttg(x, t)$ and $g(x, t)$ be the corresponding solutions to (\ref{krf}) and
(\ref{krfn}) respectively. We will show that the
eigenvalues of $Rc(p, t)$ relative to $g(p, t)$ actually converge to a fixed set
of numbers as $t\to \infty$. Here $Rc(p,t)$ is the Ricci tensor of $g(t)$ at $p$. If in addition that $(M,\wt g)$ has   maximal volume growth with positive Ricci curvature or has positive curvature operator, then we will show that for any $p \in M$, $(M, g(x, t), p)$
approaches an expanding gradient \KRS as $t\to \infty$ in the sense of limiting
solutions to the \KRF (\cite{Ha2}).
\begin{prop}\label{s3p1} Let $(M^n,g_\ijb(x))$, $\ttg(x,t)$, $g(x,t)$ be as in Theorem \ref{shitamthm} satisfying (\ref{main1e1}).  Let $p\in M$ be a fixed point in $M$ and let $\lambda_1(t)\ge
\dots>\lambda_n(t)\ge0$ be the eigenvalues of $R_\ijb(p,t)$ relative to
$g_\ijb(p,t)$. \begin{enumerate}
\item[(i)] For any $\tau>0$,
$$\frac{\det(R_\ijb(p,t)+\tau\delta_{ij})}{\det(g_\ijb(p,t))}$$ is nondecreasing in
$t$.
\item[(ii)] Assume in addition that $\tilde{g}_\ijb(x)$ has positive Ricci curvature. Then there is a constant $C>0$ such that $\lambda_n(t)\ge C$ for all
$t$.
\item[(iii)] For $1\le i\le n$ the limit $
\lim_{t\to\infty}\lambda_i(t)$ exists.
\item[(iv)] Let $\mu_1>\dots>\mu_l\ge 0$ be
the distinct limits in (iii) and let $\rho>0$ be such that
$[\mu_k-\rho,\mu_k+\rho]$, $1\le k\le l$ are disjoint. For any $t$, let $E_k(t)$
be the sum of the eigenspaces corresponding to the eigenvalues $\lambda_i(t)$
such that $\lambda_i(t)\in (\mu_k-\rho,\mu_k+\rho)$. Let $P_k(t)$ be the orthogonal projection (with
respect to $g(t)$) onto $E_k(t)$. Then there exists $T>0$ such
that if $t>T$ and if   $w\in T_p^{(1,0)}(M)$, $|P_k(t)(w)|_t$ is continuous in $t$,
where  $|\cdot|_t$
is the length measured with respect to the metric $g(p,t)$.
\end{enumerate}
\end{prop}
\begin{proof} (i): By the Li-Yau-Hamilton (LYH) inequality in \cite{cao} and in  \cite[Theorem 2.1]{Cao}, if
\begin{equation}\label{s3e3}
Z_\ijb=\frac{\p R_\ijb}{\p t}+g^{k\bar l}R_{i\bar l}R_{k\bar j}+R_\ijb
\end{equation}
then
\begin{equation}\label{s3e4}
Z_\ijb w^iw^{\bar j}\ge 0
\end{equation}
for any $w\in T^{(1,0)}(M)$. For any $\tau>0$,  denote
$$
\phi(t)= \frac{\det(R_{\ijb}+ \tau  g_{\ijb})}{\det(g_{\ijb})}
$$
at $(p,t)$. Denote $p_{\ijb}=R_{\ijb}+ \tau g_{\ijb}$ as in \cite{cao}
and note that
$(p_\ijb)$ is invertible and its inverse is denoted by $(p^{\ijb})$. We have
\begin{equation}
\begin{split}
\frac{\p}{\p t}\log \phi&= p^{\ijb}\frac{\p}{\p t}p_{\ijb}-g^{\ijb}\frac{\p}{\p t}g_{\ijb}\\
&=p^{\ijb}\lf(\frac{\p}{\p t}R_{\ijb}-\tau(R_\ijb+g_\ijb)\ri)+g^{\ijb}\lf(R_{\ijb}+g_\ijb\ri)\\
&\ge  p^{\ijb}\lf(-g^{k\bar l}R_{i\bar l}R_{k\bar j}- R_{\ijb}-\tau(R_\ijb+g_\ijb)\ri)+g^{\ijb}\lf(R_{\ijb}+g_\ijb\ri)\\
&=p^{\ijb}\lf(-g^{k\bar l}R_{i\bar l}R_{k\bar j}-(\tau+1)p_\ijb\ri)+\tau^2 p^\ijb
g_\ijb+g^\ijb(R_\ijb+g_\ijb)
\end{split}
\end{equation}
where we have used (\ref{s3e3}) and (\ref{s3e4}). Now at the point $(p,t)$, we
choose a unitary basis such that $g_{\ijb}=\delta_{ij}$ and
$R_{\ijb}=\lambda_i\delta_{ij}$. Then $p_{\ijb}=(\lambda_i+ \tau  )\delta_{ij}$
and $p^{\ijb}=(\lambda_i+\tau)^{-1}\delta_{ij}$. Hence we have
\begin{equation}
\begin{split}
\frac{\p}{\p t}\log \phi&\ge -\sum_{i=1}^n \frac{\lambda_i^2}{\lambda_i+\tau}-(\tau+1)n+\sum_{i=1}^n\frac{\tau^2}{\lambda_i+\tau}+\sum_{i=1}^n\lambda_i+n\\
&=\sum_{i=1}^n \lf(\frac{-\lambda_i^2}{\lambda_i+\tau}-\tau+\frac{\tau^2}{\lambda_i+\tau}+\lambda_i\ri)\\
&=0.\end{split}
\end{equation}
From this (i) follows.

(ii): By (i),  we conclude that $\frac{\det(R_\ijb(p,t))}{\det(g_\ijb(p,t))}$ is
nondecreasing. (This fact has been  proved in \cite{cao}.) Moreover,
$$\lim_{t\to-\infty} \frac{\det(R_\ijb(p,t)}{\det(g_\ijb(p,t))}=
\frac{\det(R_\ijb(p))}{\det(g_\ijb(p))}
$$
where the right side is in terms of the initial metric $g$ for (\ref{krf}). Since the Ricci curvature is assumed
to be positive,
 $\frac{\det(R_\ijb(p,t)}{\det(g_\ijb(p,t))}\ge C_1$ for some positive
constant $C_1$ for all $t$. On the other hand, by Corollary 2.1 there
is a constant $C_2$ independent of $t$ such that $\lambda_1(t)\le
C_2$. From these two facts, part (ii) of the proposition follows.

(iii): Choose a unitary basis $v_1,\dots,v_n$ for $T_p^{(1,0)}(M)$ with respect to
the metric $g(p,0)$. Using the Gram-Schmidt process, we can obtain a unitary
basis $v_1(t),\dots,v_n(t)$ for $g(p,t)$. Since $g(t)$ is smooth in $t$, we
conclude that the $v_i(t)$'s are smooth in $t$. That is to say,
$v_i(t)$ is a linear combination of a fixed basis of $T_p^{(1,0)}(M)$
with smooth coefficients.  Denote by $R_\ijb(t)=Rc(v_i(t),\bar
v_j(t))$ the components of $Rc(p,t)$ with respect to this basis. Then $R_\ijb(t)$ is
also smooth in $t$. By (i) and Corollary 2.1, for any $\tau>0$,
\begin{equation}\label{s3e33}
\lim_{t\to\infty}\det(R_{\ijb}(t)+\tau \delta_{ij})=c(\tau)
\end{equation}
 exists.

Now $\lambda_i(t)$ are uniformly bounded functions in $t$. To prove (iii), it is sufficient to prove that if
$t_k\to\infty$, $t_k'\to\infty$ and
$$
\lim_{k\to\infty}\lambda_i(t_k)= \tau_i,\ \lim_{k\to\infty}\lambda_i(t_k')=\tau_i'
$$
for all $i$, then $\tau_i=\tau_i'$.

By (\ref{s3e33}),  we have
$$
 \prod_{i=1}^n(\tau_i+\tau)=\prod_{i=1}^n(\tau_i'+\tau)
$$
for all $\tau>0$. Since $\tau_1\ge \dots\ge \tau_n$ and $\tau_1'\ge \dots\ge
\tau_n'$, we must have $\tau_i=\tau_i'$. This completes the proof of (iii).

(iv): By (iii), if $T$ is large enough, for each $i$ we have $\lambda_i(t)\in
(\mu_k-\rho,\mu_k+\rho)$ for some $k$ for all $t\ge T$. Hence $\dim E_k(t)$ is
constant in $t$ for $t\ge T$. Let $P_k(t)$ be the orthogonal projection (with
respect to $g(t)$) onto $E_k(t)$.  We also denote the matrix of this projection,
 with respect to the basis $v_1(t),\dots,v_n(t)$ in (iii), by $P_k(t)$. Then
$$
P_k(t)=-\frac{1}{2\pi\ii}\int_C (R_\ijb-z\delta_{ij})^{-1}dz.
$$
where $C$ is a circle on the complex plane with center at $\mu_k$ and radius
$\rho$, see \cite[p. 40]{K} for example. It is easy to see
that the matrix valued function  $P_k(t)$ is continuous in $t$. Hence (iv) is
true.
\end{proof}
\begin{rem}\label{s3r1} The facts that the scalar curvature $R(t)$ and
$\det(R_\ijb(t))/\det(g_\ijb(t))$ are nondecreasing have been proved in
\cite{cao}
\end{rem}

 Next, we will study the asymptotic behavior of the manifolds
 $(M^n,g(t))$ as $t \to \infty$. We will need the following lemma from \cite{F}:

\begin{lem}\label{uniqueness} Let $(M^n,g_\ijb)$ be a complete noncompact K\"ahler manifold with bounded curvature. Suppose there is a smooth function $f$ such that $\sqrt{-1}\p\ol\p f=Rc$. Let $g_\ijb(t)$ and $\widehat g_\ijb(t)$ be two solutions of (\ref{krf}) on $M\times [0,T]$, $T>0$ with the same initial data $g_\ijb$ such that
\begin{equation}\label{equivalent}
c^{-1}g_\ijb(x)\le g_\ijb(x,t), \widehat g_\ijb(x,t)\le cg_\ijb(x)
\end{equation}
for some constant $c>0$ for all $(x,t)\in M\times[0,T]$. Then $g_\ijb(x,t)=\widehat g_\ijb(x,t)$ on $M\times[0,T]$.
\end{lem}
In \cite{Cao} it was proved  by Cao that for any $t_k \to \infty$, if $|R(p_k, t_k)|$ is the maximum of the scalar curvature on $M$ at $t_k$, then the blow down limit of $g(t)$ along $(p_k, t_k)$ is an expanding gradient \KRS.  Recently, it is  shown by Ni in \cite{Ni0} that the result is still  true for an arbitrary sequence $p_k \in M$, $t_k\to\infty$.  In 
the special case that the sequence $p_k=p$ is fixed at an arbitrary $p\in M$, the result
follows from a rather simple observation  and the argument in \cite{Cao}, which we present below. 
\begin{prop}\label{s3p2} Assume the conditions and notation of
  Proposition \ref{s3p1}. In addition, assume the initial metric
$\wt g(x,0)= \tilde{g}_\ijb(x)$ of (\ref{krf}) has either maximal volume
  growth with positive Ricci curvature or has positive curvature
  operator.  Let $p\in M$ be a fixed point. The given any
$t_k\to\infty$, we can find a subsequence also denoted by $t_k$, a complete
noncompact complex manifold $N^n$, and a family of K\"ahler metrics $h(t)$ on $N$
satisfying (\ref{krfn}) for all $t\in \rr$ such that $(M^n,g_k(t))$,
where $g_k(t)=g(t_k+t)$ for all $t\in \rr$,  converges to $(N,
h(t))$ in the following sense:  There
exists a family of diffeomorphisms $F_k:U_k\subset N\to M$ with the following
properties.
\begin{enumerate} \item[(i)] Each $U_k$ contains $o$ where $o\in N$ is a fixed
point and $F_k(o)=p$. \item[(ii)] $U_k$ is open and the $U_k$'s exhaust $N$.
\item[(iii)] $(U_k,F_k^*(g_k(t)))$ converges in $C^\infty$ norm uniformly to
$h(t)$ in $N\times \rr$.
\end{enumerate}
Moreover $(N,h(t))$ is a gradient K\"ahler-Ricci soliton. More precisely, there
is a family of biholomorphisms $\phi_t$ of $N$ determined by the gradient of some
real valued function such that $o$ is a fixed point of each $\phi_t$  and
$\phi_t^*(h(0))=h(t)$ for all $t\ge0$.
\end{prop}

\begin{proof} The existence of $t_k$, $N$, $h(t)$ and $F_k$
satisfying (i)--(iii) is a consequence of Theorem \ref{shitamthm}  and the
compactness theorem of Hamilton \cite{Ha2}.

We now prove the last assertion in the Proposition.  Begin by noting that $\lim_{t\to\infty}R(t)$
exists by Proposition \ref{s3p1}, where $R(t)$ is the scalar curvature of $g(t)$ at $p$. Let $R^h(t)$
be the scalar curvature of $h(t)$ at $o$. Then for any $t$, $t'$
\begin{equation}\label{s3e6}
R^h(t)=\lim_{k\to\infty}R(t_k+t)=\lim_{k\to\infty}R(t_k+t')=R^h(t').
\end{equation} Now consider the metric $\wt h(t)=th(\log t)$
for $t\ge 1$ .  Then $\wt h$ is a solution to (\ref{krf}) on $N\times[0,\infty)$. Also,
since $g(t)$ has uniformly bounded   curvature in spacetime by Corollary 2.1, $h(t)$
also has uniformly bounded curvature in spacetime.  By Proposition \ref{s3p1} (ii), the Ricci
curvature of $h(t)$ at $p$ is positive. Moreover, by Theorem 2.1, the facts that
$M$ is simply connected and that the metrics $g(t)$ are decreasing in $t$, we can
conclude that $N$ is simply connected.  By \cite{cao1}, it is easy to see that
$h(t)$ and hence $\wt h(t)$ have positive Ricci curvature.  Now (\ref{s3e6})
implies that $t\ttR(t)$ is constant where $\ttR(t)$ is the scalar curvature of
$\wt h(t)$ at $p$. Hence $ \frac{\p  } {\p t}(t\ttR) =0 $ for all $t$,
and by the proof
of Theorem 4.2 in \cite{Cao}, there is a real valued function $f$ such that
$f_\ijb(x)=\ttR_\ijb(x,1)+  \wt h_\ijb(x,1)$ on $N$ with $f_{ij}\equiv0$ and
$\nabla f(o)=0$.

Let $\phi_t(x)$ be the integral curve of $-\frac12\nabla f$ on $N$ with initial
point $x$.  We claim that $\phi_t(x)$ is defined for all $x$ and $t$. Let $\wt
h_{AB}$ and $\ttR_{AB}$ be the Riemannian metric $2Re(\wt h_\ijb)$ and Ricci
curvature of $\wt h_{AB}$. Then $f_{AB}=\ttR_{AB}+\wt h_{AB}$.  Observe that
as in (\cite{Ha4} \S20), we have
\begin{equation}\label{s3e7}
|\nabla f|^2+\ttR=2f+2C_1
\end{equation}
where $\ttR$ is the scalar curvature of $\wt h(1)$ and $C_1$ is a constant.

 Now, as long as $\phi_t(x)$ is defined in  on $[-T,0]$ for $T>0$,  then for $0\le t\le T$
\begin{equation}
\begin{split}
f(\phi_{-t}(x))-f(x)&=\int_0^{-t}\frac{d}{ds}f(\phi_{s} (x))ds\\
&=\int_0^{-t}\langle \nabla f(\phi_{s}(x)),\frac{d}{ds}\phi_{s}(x)\rangle ds\\
&=\frac12\int_0^t |\nabla f(\phi_{-s}(x))|^2ds\\
&\le    \int_0^{t}  f(\phi_{-s}(x)) ds+C_1t\\
\end{split}
\end{equation}
by (\ref{s3e7}). Hence we have $f(\phi_{-t}(x))\le C_2$ for some constant
depending only on $T$, $C_1$ and $f(x)$. One can also prove that $f(\phi_t(x))\le
f(x)$ for $t>0$ as long as $\phi_t(x)$ is defined up to $t$. Since $f$ is an
exhaustion function by (\cite{CT}, Lemma 3.1), we conclude that $f(\phi_t(x))$
remains in a fixed compact set   on any bounded interval of $\rr$ as long as
$\phi_t$ is defined on that interval. From this it is easy to see that
$\phi_t(x)$ is defined for all $t$.  Since $\nabla f$ is a holomorphic
vector field, $\phi_t$ is in fact a biholomorphism on $N$ for all $t$.

Let $h_1(t)=\phi_t^*(\wt h(1))=\phi_t^*(h(0))$ and let $\wt
h_1(t)=th_1(\log t)$ for $t\ge 1$. We will show that $\wt h_1(1)=\wt h(1)$. Since $h(t)$ has nonnegative holomorphic bisectional curvature such that its scalar curvature is uniformly bounded in spacetime, $\wt h(t)$ also has nonnegative holomorphic bisectional curvature with $t\wt R(t)$ being uniformly bounded in spacetime where    $\wt R(t)$ is the scalar curvature of $\wt h(t)$.  By \cite[Theorem 2.1]{NT} and \cite[Theorem 5.1]{NST}, we can find a potential function for the Ricci tensor of $\wt h(1)$. Since the curvature of $\wt h$ and $\wt h_1$ are uniformly bounded on $M\times[0,T]$ for fixed $T>0$, it is easy to see that they satisfy (\ref{equivalent}).
By Lemma \ref{uniqueness}, we conclude that   $\wt h_1(t)=\wt h(t)$ for $t\ge 1$. Hence $h_1(t)=h(t)$ for all $t\ge 0$.  This completes the proof of the proposition.
\end{proof}

Let $t_k\to\infty$ such that $(M,g_k(t))$ converges to $(N,h(t))$ as in
Proposition \ref{s3p2}.  We will describe this convergence in terms of
the convergence of certain specific quantities.  For
simplicity, we identify $(M,g_k(t))$ near $p$ with $(U,F_k^*(g_k(t))$ for some
open set $U\subset N$ containing $o$.  Let $J_k$ be the complex structure on $U$ given by
the pullback of the complex structure of $M$ under $F_k$ and let $J$ be the complex
structure of $N$. By taking a subsequence we may also assume that $J_k\to J$. Let
$w_k\in T_p^{(1,0)}(M)$ with $|w_k|_{g_k(0)}=1$ and let
$w_k(t)=w_k/|w_k|_{g_k(t)}$ for $t\ge 0$. Denote $w_k=x_k-\ii J_k(x_k)$ where
$x_k$ is in the real tangent space of $M$ at $p$ which is identified with the
real tangent space of $N$ at $o$. Assume that $x_k\to x$. Then  $J_k(x_k)\to
J(x)$. Let $u=x-\ii J(x)$ and let $u(t)=u/|u|_{h(t)}$ for $t\ge 0$. Note that
$|u|_{h(0)}=1$.

Assume the conditions and notation of Proposition \ref{s3p2} and
Proposition \ref{s3p1}.  Then we can see
that by the propositions, the eigenvalues of the Ricci curvature of $h(t)$ with
respect to $h(t)$ at $o$ are $\mu_1>\dots>\mu_l>0$ such that the multiplicity of
$\mu_i$ is $\dim E_i(t)$ for $t$ large enough.

Let $E_i^h(t)$ be the eigenspace of the Ricci tensor of $h(t)$ corresponding to
the eigenvalue $\mu_i$.

We want to prove the following:

\begin{lem}\label{s3l1} With the assumptions as in Proposition \ref{s3p2} and with the above
notations. Suppose $w_k(t)=\sum_{i=1}^l w_{k,i}(t)$ where $w_{k,i}(t)$ is the
orthogonal projection of $w_k(t)$ onto $E_{i}(t+t_k)$ with respect to
$g_k(t)=g(t_k+t)$ and suppose  $u(t)=\sum_{i=1}^l u_i(t)$ where $u_i(t)$ is the
orthogonal projection of $u(t)$ onto $E_i^h(t)$ with respect to $h(t)$. Then for any $T>0$,  the
following are true:
\begin{enumerate} \item[(i)] $w_k(t)$ converges uniformly to
$u(t)$ on $t\in [0,T]$ in the sense that the real parts and the imaginary parts
of $w_k(t)$ converge uniformly to the real part and imaginary part of $u(t)$
respectively.
\item[(ii)] $Rc^k_t(w_k(t),\bar w_k(t))$ converges uniformly to
$Rc^h_t(u(t),\bar u(t))$ on $t\in [0,T]$ where $Rc_t^k$ is the Ricci tensor of $g_k(t)$ at $p$
and $Rc_t^h$ is the Ricci tensor of $h(t)$ at $o$.
\item[(iii)]  By passing to a
subsequence if necessary, for $1\le i\le l$,  $ |w_{k,i}(t)|_{g_k(t)}$  converge
uniformly to $|u_i(t)|_{h(t)}$ on $t\in [0,T]$.
\end{enumerate}
\end{lem}
\begin{proof} (i): Since $g_k(t)$ converges uniformly to $h(t)$ on $[0,T]$ at $o$
and since $w_k\to u$, $|w_k|_{g_k(t)}$ converge to $|u|_{h(t)}$ uniformly on
$[0,T]$. From this it is easy to see that (i)  is true.

(ii): Since $g_k(t)$ converges uniformly on $U\times[0,T]$ in $C^\infty$ norm, by
(i) it is easy to see that (ii) is true.

(iii): Let $v_{k}^{(1)},\dots,v_{k}^{(n)}$ be a unitary basis for $T_p^{(1,0)}(M)$ with
respect to $g_k(0)$. Passing to a subsequence if necessary, we may assume that
they converge to a unitary basis $u^{(1)} ,\dots, u^{(n)}$ of $T_o^{(1,0)}(N)$ with
respect to $h(0)$. Using the Gram-Schmidt process, we can obtain
$v_{k}^{(1)}(t),\dots,v_{k}^{(n)}(t)$ to be a unitary basis for $T_p^{(1,0)}(M)$ with
respect to $g_k(0)$ and   a unitary basis $u^{(1)}(t),\dots,
u^{(n)}(t)$ of $T_o^{(1,0)}(N)$ with respect to $h(t)$.  We claim that
$v_{k}^{(i)}(t)$ converges
to $u^{(i)}(t)$ uniformly on $[0,T]$. In fact, since $g_k(t)$ converge to $h(t)$ uniformly on $[0,T]$ and $v_k^{(1)}\to u^{(1)}$, $|v_k^{(1)}|_{g_k(t)}\to |u^{(1)}|_{h(t)}$ uniformly on $[0,T]$. Define $ v_k^{(1)}(t)=v_k^{(1)}/|v_k^{(1)}|_{g_k(t)}$ and $u^{(1)}(t)=u^{(1)}/|u^{(1)}|_{h(t)}$. Then $v_{k}^{(1)}(t)$   converge
to $u^{(1)}(t)$ uniformly on $[0,T]$. Suppose we have found $v_{k}^{(i)}(t)$, $1\le i\le m$ and $u^{(i)}(t)$, $1\le i\le m$ such that (a) $v_{k}^{(i)}(t)$, $1\le i\le m$ are unitary with respect to $g_k(t)$  and are linear combinations  of  $v_{k}^{(i)} $, $1\le i\le m$; (b) $u^{(i)}(t)$, $1\le i\le m$ are unitary with respect to $h(t)$ and are linear combinations of
$u^{(i)} $, $1\le i\le m$; and (c) $v_{k}^{(i)}(t)$   converge
to $u^{(i)}(t)$ uniformly on $[0,T]$ for $1\le i\le m$. Define
$$
v_{k}^{(m+1)}(t)= \frac{v_{k}^{(m+1)} -\sum_{i=1}^m \langle v_{k}^{(m+1)},  v_{k}^{(i)} (t)\rangle_{g_k(t)}v_{k}^{(i)} (t)}{|v_{k}^{(m+1)} -\sum_{i=1}^m \langle v_{k}^{(m+1)},  v_{k}^{(i)} (t)\rangle_{g_k(t)}v_{k}^{(i)} (t)|_{g_k(t)}}
$$
and define
$$
u^{(m+1)}(t)= \frac{u^{(m+1)} -\sum_{i=1}^m \langle u^{(m+1)},  u^{(i)} (t)\rangle_{h(t)}u^{(i)} (t)}{|u^{(m+1)} -\sum_{i=1}^m \langle u^{(m+1)},  u^{(i)} (t)\rangle_{h(t)}u^{(i)} (t)|_{h(t)}}
$$
Then it is easy (a), (b) and (c) are still true with $m$ replaced by $m+1$. Hence by induction, we can construct $v_k^{(i)}(t)$ and $u^{(i)}(t)$ as claimed.

Let
$R^k_{\ijb}(t)=Rc^k_t(v_{k}^{(i)}(t),\bar v_{k}^{(j)}(t))$ and let
$R^h_\ijb(t)=Rc^h_t(u^{(i)}(t),\bar u^{(j)}(t))$. Then as in (ii), we can prove that
$R^k_{\ijb}(t)$ converge to $R^h_\ijb(t)$ uniformly on $[0,T]$. Denote $P^k_i(t)$
to be the matrix with respect to the basis $v_{k}^{(1)}(t),\dots,$ $v_{k}^{(n)}(t)$ of the
orthogonal projection onto $E_{i}(t+t_k)$ with respect to $g_k(t)$. Denote
$P_i(t)$ to be the matrix with respect to the basis $u^{(1)}(t),\dots,u^{(n)}(t)$ of the
orthogonal projection onto $E_i^h(t)$ with respect to $h(t)$. As in the proof of
Proposition \ref{s3p1}(iv),

\begin{equation}\label{s3e8}
P^k_s(t)=-\frac{1}{2\pi\ii}\int_C (R^k_\ijb(t)-z\delta_{ij})^{-1}dz
\end{equation}
and
\begin{equation}\label{s3e9}
P_s(t)=-\frac{1}{2\pi\ii}\int_C (R^h_\ijb(t)-z\delta_{ij})^{-1}dz
\end{equation}
 where $C$ is a circle on the complex plane with center at $\mu_s$ and radius $\rho$.
  Since $R^k_\ijb(t)$ converge to $R^h_\ijb(t)$ uniformly on $[0,T]$,
(iii) follows from (\ref{s3e8}),  (\ref{s3e9}) and (i).

\end{proof}

\section{Asymptotic behavior of K\"ahler Ricci flow (II)}
Let $(M^n,\ttg)$ be as in Theorem 2.1 with maximal volume growth or with positive curvature operator and let $g(x,
t)$ be the corresponding solution to (\ref{krfn}).  As before, we denote the
eigenvalues of $Rc(p,t)$ by $\lambda_i (t)$ for $i=1,...,n$ and we let $\mu_k$,
$E_k(t)$ and $P_k(t)$ for $k=1,...,l$ be as in Proposition \ref{s3p1}. We let
$n_m$  for $m=0,...,l-1$ be such that $\lambda_k (t) \in (\mu_m -\rho,
\mu_m+\rho)$ for all $n_{m}\leq k\leq n_{m+1}$ and $t$ sufficiently large such
that the intervals $[\mu_m
  -\rho, \mu_m+\rho]$ are disjoint as in Proposition \ref{s3p1} part (iv).  For any nonzero vector $v \in T_p ^{1, 0}
(M)$, let $v(t)=v/ |v|_t $ where $|v|_t$ is the length of $v$ with respect to $g(t)$ and $v_i(t)=P_i(t)v(t)$.

The goal of this section will be to prove that $Rc(p,t)$ can be
`diagonalized' simultaneously near infinity in a certain sense and
that $g(t)$ is `Lyapunov regular', to borrow a notion from dynamical
systems (see \cite{BP}).

In the following  lemmas we assume  that the initial
metric $\wt g(0)$ in (\ref{krf}), and thus by Proposition \ref{s3p1} $g(x,t)$ for
all $(x, t)$, has positive Ricci curvature.

Now let $(N, h(t))$ be a gradient \KRS as in Proposition \ref{s3p2} and
let $o\in N$, $\phi_t$ and $E_i^h (t)$ also be as in the Proposition.
For any nonzero vector $w \in T_o ^{1, 0}(N)$ let $w(t)=w/ |w|_{h(t)}$ and
$w_i(t)$ be the projection of $w(t)$ onto $E_i^h (t)$.  We begin by
making the following observation.

Let $\phi_t$ be given by   $-\frac12\nabla f$ such that $f_{\ijb}(x)=R^h_\ijb(x,0)+h_\ijb(x,0)$ and $f_{ij}=0$. Near $o$, we may choose local coordinates  $z_i$ such that $\p_i=\frac{\p}{\p z_i}$ are unitary at $o$ which diagonalize  $f_\ijb$ at $o$. We also assume that the origin corresponds to $o$.  Then $\mu_1>\mu_2>\dots>\mu_l>0$ are distinct eigenvalues of $Ric^h$ at $t=0$ with respect to $h(0)$.   Since $\p_i$ are  eigenvectors of $f_\ijb$, for each $i$ we have
\begin{equation}\label{basis}
(\phi_t)_*(\p_i)=e^{-\frac12(\mu_j+1)t}\p_i
\end{equation}
for some $j$ at $o$.
 Because of (\ref{basis}) and the fact that $\p_i$ are also eigenvectors of $R_\ijb$ at $o$ and $t=0$, $E_i^h(0)=E_i^h(t)$
and   $w_i(t)= w_i(0)/|w|_{h(t)}$.

\begin{lem}\label{ns4l1} Let $(N, h(t))$ be a gradient K\"ahler-Ricci soliton and $w\in T_o^{(1,0)}(N)$ with $|w|_{h(0)}=1$ as
  above. Let $1\le m<l$, and suppose $a<\sum_{j=m+1}^l|w_j(0)|_{h(0)}^2<1-a$ for some $0<a<1$.  Then for $t\ge 0$,
$$
\frac{\sum_{j=m+1}^l|w_j(t)|_{h(t)}^2}{\sum_{j=1}^m|w_j(t)|_{h(t)}^2}\ge
\frac{\sum_{j=m+1}^l|w_j(0)|_{h(0)}^2}{\sum_{j=1}^m|w_j(0)|_{h(0)}^2}\cdot
e^{(\mu_m-\mu_{m+1})t}.
$$
In particular,
$$
\sum_{j=m+1}^l|w_j(t)|_{h(t)}^2\ge \sum_{j=m+1}^l|w_j(0)|_{h(0)}^2
$$
for $t\ge 0$.  Moreover, for any $\delta>0$, there is a $t_0$ depending only on the $a$, $\mu_m$, $\mu_{m+1}$ and $\delta$ such that for all $t\ge t_0$,
$$
 \sum_{j=m+1}^l|w_j(t)|_{h(t)}^2\ge 1-\delta.
$$\end{lem}
\begin{proof} For simplicity, let us denote $|\cdot|_{h(t)}$ simply by $|\cdot|_t$.

\begin{equation}
\begin{split}
|w_j(t)|_t^2 &= \frac{|(\phi_t)_*(w_j(0))|^2_0}{|(\phi_t)_*(w)|_0^2}\\
&=\frac{e^{(- \mu_j-1)t}|w_j(0)|^2_0}{|w|_t^2}.
\end{split}
\end{equation}
Hence for $t\ge0$
\begin{equation}\label{nos4e1}
\sum_{j=1}^m|w_j(t)|_t^2=\frac{\sum_{j=1}^{m}e^{(-
\mu_j-1)t}|w_j(0)|^2_0}{|w|_t^2}\le \frac{e^{(-
\mu_m-1)t}\sum_{j=1}^{m}|w_j(0)|^2_0}{|w|_t^2}
\end{equation}
because $\mu_1>\dots>\mu_l$. Similarly,
\begin{equation}\label{nos4e2}
\sum_{j=m+1}^l|w_j(t)|_t^2=\frac{\sum_{j=m+1}^{l}e^{(-
\mu_j-1)t}|w_j(0)|^2_0}{|w|_t^2}\ge \frac{e^{(-
\mu_{m+1}-1)t}\sum_{j=m+1}^{l}|w_j(0)|^2_0}{|w|_t^2}.
\end{equation}
The lemma then follows from (\ref{nos4e1}) and (\ref{nos4e2}).
 \end{proof}

Because of Proposition \ref{s3p2} and Lemma \ref{s3l1}, we expect to have similar behavior for $g(t)$ for $t$ large. More precisely, we have the following:
\begin{lem}\label{ns4l2}
Let $v_k \in T_p^{(1, 0)} (M)$ be a sequence such that $|v_k|_0 =1$ for each $k$.  Let
$t_k \to \infty$ be a sequence in time.  Define $f_{ik}(t):=|P_i(t)v_k (t)|_t^2$.
\begin{enumerate}
\item[(i)] Suppose there exists $a >0$ and $1\le m\le l$ for which
\begin{equation}
\sum_{i\geq m} f_{ik}(t_k) =a
\end{equation}
for all $k$.  Then for any sequence $s_k > t_k$ we have
\begin{equation}
\liminf_{k\to \infty} \sum_{i\geq m} f_{ik}(s_k) \geq a.
\end{equation}
\item[(ii)] Suppose there exists $1>a >0$ and $1\le m\le l$ for which
\begin{equation}\label{ns4e.01}
a \leq \sum_{i\geq m} f_{ik}(t_k) \leq 1-a.
\end{equation}
for all $k$.  Then for any $1>\delta>0$ there exists $T>0$ such that
\begin{equation}\label{ns4e.02}
\liminf_{k \to \infty} \sum_{i\geq m}f_{ik}(t_k +T) \geq 1-\delta.
\end{equation}
\end{enumerate}
\end{lem}

\begin{proof}
Suppose (i) is false. Then $m>1$ and there exists a subsequence of $t_k$ which we will
also denote by $t_k$, a sequence $s_k >t_k$, and some $\e>0$ for which
\begin{equation}\label{ns4e0}
\sum_{i\geq m} f_{ik}(s_k) \leq a-\e
\end{equation}
for all $k$.  Thus by the continuity of $f_{ik} (t)$ in $t$ for each $i$ (see Proposition \ref{s3p1}(iv)), there is a sequence
$t_k<T_k<s_k$ such that
\begin{equation}\label{ns4e1}
\sum_{i\geq m} f_{ik}(T_k)=a-\frac{\e}{2}
\end{equation}
and
\begin{equation}\label{ns4e2}
\sum_{i\geq m} f_{ik}(t) \leq a-\frac{\e}{2}
\end{equation}
for all $t\in [T_k, s_k]$.

Now define $g_k(t)=g(T_k+t)$. Then we may assume that $(M,g_k(t))$ converges to a
soliton $(N,h(t))$ as in Proposition \ref{s3p2} such that $p$ corresponds to the
stationary point $o$. We may also assume that $v_k(T_k)$ converges to a vector $w$ in
$T^{1, 0}_o (N)$ where $w$ has length 1 in with respect to $h(0)$.  Then by Lemma
\ref{s3l1}(iii), for any $T>0$,  we have
\begin{equation}\label{ns4e3}
\lim_{k\to \infty} \sum_{i\geq m} f_{ik}(T_k +t)=\sum_{i\geq m}|w_i(t)|_{h(t)}
\end{equation}
uniformly for all $t\in [0, T]$, where $w(t)=w/|w|_{h(t)}$ and $w_i(t)$ is the orthogonal projection of $w(t)$ onto the eigenspace of $Ric^h(t)$ at $o$ of the eigenvalue $\mu_i$ with respect to $h(t)$.

We claim that  $s_k -T_k > \tau$ for
some $\tau >0$. Otherwise, we may assume that $s_k -T_k \to 0$, and
thus from (\ref{ns4e0}), (\ref{ns4e1}) and (\ref{ns4e3}) we may draw
the contradiction that
$$
a-\frac\e2=\sum_{i\ge m}|w_i(0)|_{h(0)}\le a-\e.
$$
This proves the claim.
Thus from (\ref{ns4e1}), (\ref{ns4e2}) and (\ref{ns4e3}) we may conclude
that
\begin{equation}\label{ns4e4}
\sum_{i\geq m} w_i(0)=a-\frac\e2
\end{equation}
and
\begin{equation}\label{ns4e5}
\sum_{i\geq m} w_i(t) \leq a-\frac{\e}{2}
\end{equation}
for all $t\in [0, \tau]$.  But (\ref{ns4e4}) and (\ref{ns4e5}) contradict Lemma
\ref{ns4l1}. This completes the proof of (i) by contradiction.

We now suppose (ii) is false. Note that $m>1$  because $0<a<1$. Then there exists a $\delta>0$ with the property
that: given any $T>0$, there exists a subsequence of $t_k$, which we also denote
by $t_k$, for which
\begin{equation}\label{ns4e6}
\sum_{i\geq m}f_{ik}(t_k +T) \leq 1-\delta.
\end{equation}
for all $k$.

 Now we define $g_k(t)=g(t_k+t)$ and assume $(M,g_k(t))$ converges to a
soliton $(N,h(t))$ as in the proof of  (i).  We also assume that $v_k(t_k)$ converges to a vector $w$
in $T^{1, 0}_o( N)$ where $w$ has length 1 with respect to $h(0)$. Then by taking a limit
as in the proof of  (i), using Lemma \ref{s3l1}(iii), (\ref{ns4e.01}) and (\ref{ns4e6}), we have
\begin{equation}\label{ns4e7}
a \leq \sum_{i\geq m} w_i(0) \leq 1-a
\end{equation}
and
\begin{equation}\label{ns4e8}
\sum_{i\geq m}w_i(T) \leq 1-\delta.
\end{equation}
But for $T$ sufficiently large depending only on $a$, $\mu_{m-1}$, $\mu_m$ and  $\delta$, (\ref{ns4e7}) and
(\ref{ns4e8}) contradict Lemma
\ref{ns4l1}.  This complete our proof of (ii) by contradiction.
\end{proof}

We are ready to prove the main theorem in this
section. \begin{thm}\label{ns4t1} Let $(M^n,\ttg)$ be as in Theorem
  2.1 with maximal volume growth or with positive curvature operator
  and let $g(x,t)$ be the corresponding solution to (\ref{krfn}). With
the same notation as above,
   $V=T_p^{(1,0)}(M)$ can be decomposed orthogonally with respect to $g(0)$ as
 $V_1\oplus \cdot\cdot\cdot\oplus V_l$ so that the following are true:
\begin{enumerate}
\item[(i)] If $v$  is a nonzero vector in  $V_i$ for some $1\le i\le l$, then
 $\lim_{t \to \infty} |v_i(t)|=1$ and thus $\lim_{t\to\infty}Rc(v(t),\bar v(t))=\mu_i$ and
$$\lim_{t \to \infty} \frac{1}{t}\log \frac{|v|_t^2}{|v|_0^2}=-\mu_i-1.$$  Moreover, the convergences are uniform over all
 $v\in V_i\setminus\{0\}$.
\item[(ii)] For $1\le i, j\le l$ and for nonzero vectors  $v \in V_i$ and $w \in V_j$ where $i\neq j$, $\lim_{t\to
\infty}\langle v(t),w(t)\rangle_t=0$ and the convergence is uniform over all such nonzero vectors $v, w$.
\item[(iii)] $\dim_\C(V_i)=n_{i}-n_{i-1}$  for each $i$.
\item[(iv)] $$\sum_{i=1}^l(-\mu_i-1)\dim_\C
V_i=\lim_{t\to\infty} \frac1t\log \frac{\det(g_{i\bar j}(t))}{\det ({g}_{i\bar
j}(0)}.$$ \end{enumerate}
\end{thm}

\begin{proof}  We first assume that the
initial metric  $\wt g(0)$ in (\ref{krf}), and thus
  $g(x,t)$ for all $(x, t)$, has positive Ricci curvature by Proposition \ref{s3p1}.

To prove (i),  let $v\in T_p (M)$ be a fixed nonzero vector and let
    $f_i(t)=|v_i (t)|^2$.  We claim that $\lim_{t \to \infty}
    f_m(t)=1$ for some $m$, and thus $\lim_{t \to \infty} f_k(t)=0$
    for all $k\neq m$.  To prove our claim it will be sufficient to
    prove the following for every $m$: Suppose $\lim_{t\to \infty}
    f_j(t)=0$ for all $j<m$.  Then either
\begin{equation}\label{ns4e10}
\lim_{t \to \infty} f_m(t)=1
\end{equation}
or
\begin{equation}\label{ns4e11}
\lim_{t \to \infty} f_m(t)=0.
\end{equation}
If $m=l$, then we must have $\lim_{t \to \infty} f_m(t)=1$. Suppose $1\le m<l$ and  $\lim_{t\to \infty}f_j(t)=0$ for all $j<m$ and that neither
  (\ref{ns4e10}) nor (\ref{ns4e11}) holds.  By the continuity of $f_i(t)$,
we can find $t_k\to\infty$ such that
\begin{equation}\label{c1}
a\le \sum_{i\ge m+1}f_i(t_k)\le 1-a
\end{equation} for some $0<a<1$.
 By letting   $v_k=v$
  for all $k$, it follows from Lemma \ref{ns4l2}(ii), we can find $T>0$,  such that passing to a subsequence if necessary we have
\begin{equation}\label{c2}
\sum_{i\ge m+1}f_i(t_k+T)\ge 1-\frac a2.
\end{equation}
For each $j$, we can find $k_j$ such that $t_{k_j}>t_j+T$. Since
$$ \sum_{i\ge m+1}f_i(t_j+T)\ge 1-\frac a2
$$
and
$$
 \sum_{i\ge m+1}f_i(t_{k_j})\le 1-a
$$
for all $j$, we may derive a contradiction from part (i) of Lemma \ref{ns4l1}.  Thus our initial assumption was false, and for any
  $v\in T_p (M)$ and $m$, either (\ref{ns4e10}) or (\ref{ns4e11})
  holds.  Thus for any nonzero $v\in T_p (M)$ we have $\lim_{t \to \infty}f_m(t)=1$ for some $m$

   Now suppose $\lim_{t \to \infty}f_m(t)=1$.  Using (\ref{krfn}), Proposition \ref{s3p1}, the definition of $\mu_i$ and the definition of $f_i(t)$,
a straight forward calculation gives
$$
\lim_{t\to\infty}\frac1t \log|v|^2_t=-\mu_m-1
$$
 Note that if
$$
\lim_{t\to\infty}\frac1t \log|v|^2_t=-\mu_i-1
$$
and
$$
\lim_{t\to\infty}\frac1t \log|w|^2_t=-\mu_j-1
$$
and $i\le j$ (so that $-\mu_j\ge-\mu_i$), then
\begin{equation}\label{s4e.1}
\lim_{t\to\infty}\frac1t \log|av+bw|^2_t\le-\mu_j-1.
\end{equation}
provided   $av+bw\neq0$.

Let $V_1$ be the subspace of $V=T_p^{(1,0)}(M)$ defined by
$$
V_1=\{v\in V\setminus\{0\}|\ \lim_{t\to\infty}\frac1t \log|v|^2_t=-\mu_1-1\}\cup\{0\}.
$$
It is easy to see that $V_1$ is a subspace by (\ref{s4e.1}). Let $V_1^\perp$ be the
orthogonal complement of $V_1$ with respect to $g(0)$. Then by the definition of $V_1$,  for any nonzero $v\in
V_1^\perp$, we have
$$
\lim_{t\to\infty}\frac1t \log|v|^2_t=-\mu_j-1
$$
for some $j>1$. Define
$$
V_2=\{v\in V_1^\perp \setminus\{0\}|\ \lim_{t\to\infty}\frac1t \log|v|^2_t=-\mu_2-1\}\cup\{0\}.
$$
Continuing in this way, we can decompose $V$ as $V=V_1\oplus \dots\oplus V_l$
orthogonally with respect to  $g(0)$, such that if $v\in V_m$, then
\begin{equation}\label{riccilimit}
\lim_{t\to\infty}f_m(t)=1,
\end{equation}
and
\begin{equation}\label{metriclimit}
\lim_{t\to\infty}\frac1t \log\frac{|v|^2_t}{|v|^2_0}=-\mu_m-1.
\end{equation}
It remains to prove that both  convergences are uniform on $V_m\setminus\{0\}$.  It is sufficient to prove the convergence in (\ref{riccilimit}) is uniform. Suppose the convergence is not uniform over $V_m\setminus\{0\}$. Then there exist $v_k\in V_m$, $t_k\to\infty$, $\e>0$ such that $|v_k|_0=1$, $v_k$ converge to some vector $v\in V_m$ and
\begin{equation}\label{uniform1}
f_{mk}(t_k)=|P_m(t_k)v_k(t_k)|^2_{t_k}\le 1-5\e.
\end{equation}
Since $f_{mk}(t)=|P_m(t)v_k(t)|_t^2\to 1$ as $t\to\infty$ for all $k$, we can find $r_k>t_k$ such that
\begin{equation}\label{uniform2}
f_{mk}(r_k)\ge  1-\e.
\end{equation}
On the other hand, for each fixed $s$, $\lim_{k\to\infty}
f_{mk}(s)=|P_m(s)v(s)|^2_s$. Moreover,
$\lim_{s\to\infty}|P_m(s)v(s)|^2_s=1$ because $v\in V_m$ and
$|v|_0=1$. Hence passing to a subsequence if necessary, we can find
$s_k\to\infty$ such that $s_k<t_k$ and
\begin{equation}\label{uniform3}
f_{mk}(s_k)\ge  1-\e.
\end{equation}
Now we claim that there exists $k_0$ such that if   $k\ge k_0$ then
\begin{equation}\label{uniform4}
\sum_{i\ge m}f_{ik}(t)\ge  1-2\e
\end{equation}
for all $t>s_k$. Otherwise, we can find $s_k'>s_k$ for infinitely many $k$ such that
\begin{equation}\label{uniform5}
\sum_{i\ge m}f_{ik}(s_k')\le  1-2\e.
\end{equation}
But (\ref{uniform3}), (\ref{uniform5}) and the fact that $s_k'>s_k$ contradicts  Lemma \ref{ns4l2}(i).

If $m=l$, then for $k\ge k_0$ (\ref{uniform4}) contradicts (\ref{uniform1}) because $t_k>s_k$. Suppose $1\le m<l$, then by (\ref{uniform1}) and (\ref{uniform4}), for $k\ge k_0$,  we have
\begin{equation}\label{uniform6}
\sum_{i\ge m+1}f_{mk}(t_k)\ge  3\e.
\end{equation}
and from (\ref{uniform2})
\begin{equation}\label{uniform7}
\sum_{i\ge m+1}f_{mk}(r_k)\le   \e.
\end{equation}
for $k$ large enough. Since $r_k>t_k$, (\ref{uniform6}) and (\ref{uniform7}) contradicts Lemma \ref{ns4l2}(i) again. This completes  the proof of part (i).

Part (ii) of the theorem  follows directly from the definition of $v(t)$
and $w(t)$, the orthogonality of the spaces $E_i(t)$ with respect to $g(t)$ and
part (i).

To prove (iii),   we begin by showing the following: Fix $1\le m\le l$.  Let
$v_k\in E_1(s_k)+\dots+E_{m}(s_k)$ with $s_k\to\infty$ such that $|v_k|_0=1$ and $v_k$
converge  to a vector $u\in T_p^{1,0}(M)$ of unit length with respect to $g(0)$.
 Then
\begin{equation}\label{ns4e12}
\lim_{t\to\infty}|u_j(t)|_t=0
\end{equation}
 for all $j>m$, where $u_j(t)=P_j(t)u(t)$ and $u(t)=u/|u|_t$ as before.

 Suppose this is false.  Then by  (i), we have
\begin{equation}
\lim_{t\to\infty}\sum_{j\ge m+1}|u_j(t)|^2_t=1.
\end{equation}
   Let $f_{jk}(t)=|P_j(t) v_k(t)|_t^2$. Since for fixed $t$,
$$
\lim_{k\to\infty} f_{jk}(t)=|u_j(t)|_t^2,
$$
as before,   given any $\frac12>\epsilon>0$ we may choose a
subsequence of $s_k$ also denoted by $s_k$, and a
 sequence $t_k<s_k$ for which $t_k \to \infty$ and
\begin{equation}\label{ns4e13}
\sum_{j\ge m+1}f_{jk}(t_k)\geq 1-\epsilon.
\end{equation}
for all $k$.  But $\sum_{j\ge m+1}f_{jk}(s_k)=0$ for all $k$ by definition.     This is impossible by   Lemma \ref{ns4l2}(i).  Thus
(\ref{ns4e12}) is true for all $j>m$.

 We now show that   for all $1\le m\le l$, $\dim_\C V_m=n_{m}-n_{m-1}$ which is equal to $\dim_\C E_m(t)$ for $t$ large enough.
Let $d_i=\dim V_i$. We claim that for any $1\le m\le l$,
\begin{equation}\label{s4e6}
d_1+\dots+d_m\ge n_m.
\end{equation}
Fix $1\le m\le l$. Choose $t_k\to\infty$. We may assume that  $\dim
E_j(t_k)=n_j-n_{j-1}$ for all $j$ and $k$. Hence  we can choose a basis  $v_1(t_k),\dots,v_{n_m}(t_k)$ of $\sum_{j=1}^mE_j(t_k)$. Using Gram-Schmidt process, we may assume
that $$ v_{1}(t_k)/|v_{1}(t_k)|_{g(0)},\dots,v_{n_m}
(t_k)/|v_{n_m}(t_k)|_{g(0)}$$ are unitary with respect to $g(0)$. Moreover, we may assume
that for $k\to\infty$, $v_j (t_k)/|v_j(t_k)|_0$ converge to some $w_j$ for all
$1\le j\le n_m$.  Hence we have $n_m$ vectors $w_{1},\dots,w_{n_m}$ They satisfy
the following:
\begin{enumerate} \item[(a)] They are unitary with respect to $g(0)$ by construction.
\item[(b)] For each $1\le j\le n_m$
$$
\lim_{t\to\infty}\frac1t\log|w_j(t)|^2_t\le -\mu_m-1
$$
by (\ref{ns4e12}).
\end{enumerate}
For each $w_j$ ($1\le j\le n_m$), $w_j=\sum_{k=1}^{l}w_{j,k}$ where $w_{j,k}\in
V_k$. If there is a $k>m$ such that $w_{j,k}\neq0$, then by (i) and the fact that $-\mu_k>-u_m$ and the definition of $V_k$,   we  have
$$
\lim_{t\to\infty}\frac1t\log|w_j(t)|^2_t\ge  -\mu_k-1>-\mu_m-1,
$$
contradicting (a).  Thus $w_j\in V_1\oplus\dots\oplus V_m$ for $1\le j\le n_m$.
From this (\ref{s4e6}) follows because the $w_j$ are linearly independent by (a).

Choose a unitary basis $v_{j,1},\dots, v_{j,d_j}$ of $V_j$ with respect to $g(0)$ for all
$1\le j\le l$. This gives a unitary basis with respect to $g(0)$ for $T_p^{(1,0)}(M)$. Let
$g_{i\jbar}(t)$ be components of $g(t)$   with respect to this basis. Then
$$
\det(g_{\ijb}(t))\le \prod_{j=1}^l\prod_{k=1}^{d_j}|v_{j,k}|_{g(t)}^2.
$$
Since $\lim_{t\to\infty}R(t)=\sum_{j=1}^l(n_j-n_{j-1}) \mu_j$ by Proposition \ref{s3p1} where $R(t)$ is the scalar curvature of $g(t)$, by (\ref{krfn}) and the above inequality we have
\begin{equation}\label{s4e7}
\begin{split}
\sum_{j=1}^l (n_j-n_{j-1})(-\mu_j-1)&=-\lim_{t\to\infty}R(t)-n\\
&= \lim_{t\to\infty}\frac1t\log \frac{\det(g_{\ijb}(t))}{\det(g_{\ijb}(0))}\\
&\le \sum_{j=1}^l\sum_{k=1}^{d_j}\lim_{t\to\infty}\frac1t\log |v_{j,k}|_{g(t)}^2\\
&=\sum_{j=1}^l d_j(-\mu_j-1).
\end{split}
\end{equation}
Let us denote $n_j-n_{j-1}$ by $k_j$, then we have
$$
\sum_{j=1}^l k_j(-\mu_j)\le \sum_{j=1}^l d_j(-\mu_j)
$$
and $\sum_{j=1}^m d_j\ge \sum_{j=1}^mk_j$ for all $1\le m\le l$ by (\ref{s4e6}). Also
$\sum_{j=1}^ld_j=\sum_{j}^l k_j=n$. Since $-\mu_1<-\mu_2<\dots<-\mu_l$, we must
have $d_j=k_j$ for all $j$. In fact, if this is not the case, since $d_1\ge k_1$, and $\sum_{j=1}^m d_j\ge \sum_{j=1}^mk_j$ for all $1\le m\le l$,  then we can find  $m$ to be the first $m$ such that $d_m>k_m$ and $d_j=k_j$ for $j<m$.
We have
\begin{equation}\label{s4e8}
\begin{split}
\sum_{j=1}^ld_j(-\mu_j)&=\sum_{j<m}k_j(-\mu_j)+k_m(-\mu_m)+(d_m-k_m)(-\mu_m)+\sum_{j>m}d_j(-\mu_j)\\
&<\sum_{j\le m}k_j(-\mu_j)+(d_m-k_m+d_{m+1})(-\mu_{m+1})+\sum_{j>m+1}d_j(-\mu_j)
\end{split}
\end{equation}
because $-\mu_m<-\mu_{m+1}$ and  $d_m-k_m>0$. If we let $d_j'=k_j$ for $1\le j\le
m$, $d_j'=d_j$ for $j>m+1$, and $d_{m+1}'=d_m-k_m+d_{m+1}$ then we have
$$
\sum_{j=1}^l k_j(-\mu_j)< \sum_{j=1}^l d_j'(-\mu_j)
$$
and $\sum_{j=1}^p d_j'\ge \sum_{j=1}^pk_j$ for all $1\le p\le l$ by (\ref{s4e6}). Also
$\sum_{j=1}^ld_j'=\sum_{j}^l k_j=n$. Moreover, $d_j'=k_j$ for all $1\le j\le m$.
By induction, we will end up with
$$
\sum_{j=1}^l k_j(-\mu_j)< \sum_{j=1}^l k_j (-\mu_j)
$$
which is impossible. This completes the proof of part (iii).

Part (iv)   follows directly from part  (iii) and the first two equalities in (\ref{s4e7}).

We have thus proved that Theorem in the case that $(M, g)$ satisfied
the additional assumption of positive Ricci curvature.  Now if the Ricci curvature is not strictly positive on $M$, we can use the
results in \cite{cao1} to reduce back to the case of positive Ricci curvature.
This completes the proof of the theorem is.
\end{proof}

\section{Uniformization}

Let $(M^n,\wt g)$ be as in Theorem \ref{shitamthm} and assume $(M,\wt g)$ has
either maximum volume growth or positive curvature operator. Let
$\ttg(t)$ be the solution of the K\"ahler-Ricci flow (2.1) and let $g(t)$ be the
corresponding solution of the normalized flow (2.3).  Fix a
point $p\in M$.  Then by Corollary \ref{Assumption}, there exist $1>r_1, r_2>0$ such that for all
$t>0$, there is a holomorphic map $\Phi_t:D(r_1)\to M$ (where $D(r_1)=\{z\in
\C^n|\ |z|<r_1\}$), satisfying the following:
\begin{equation}\label{s5e1}
\begin{cases}&\text{$\Phi_t$ is biholomorphism from $D(r_1)$ onto its image.}\\
& \Phi_t(0)=p.\\
&\text {$\Phi_t^*(g(t))(0)=g_\e$, where $g_\e$ is the standard Euclidean metric of $\C^n$.}\\
&\text{$\frac1{r_2}g_\e\le\Phi_t^*(g(t))\le r_2 g_\e$ in $D(r_1)$.}\end{cases}
\end{equation}

Let $T>0$ and let $F_{i+1}=\Phi_{(i+1)T}^{-1}\circ \Phi_{iT}$.  Then
for each $i$,
$F_i$ is a
holomorphic map from $D(r_1)$ into $\C^n$ and is a biholomorphism onto its image.
 Let $A_i=F_i'(0)$ be the Jacobian matrix of $F_i$ at $0$.   By Corollary \ref{Assumption}, we will choose $T>0$   large enough such that
\begin{equation}\label{s5e2}
F_i(D(r_1))\subset D(r_1), \ |F_i(z)|\le \delta |z|\text{\ for some
$0<\delta<1$}.
\end{equation}
Since $R_\ijb\ge0$ for all $t$, we have
\begin{equation}\label{s5e3}
a|v|\le |A_i(v)|\le b|v|
\end{equation}
for some $0<a<b<1$ for all $i$. Here $a, b, \delta$ are independent of $i$.
We will now modify(decompose) the maps $F_i$ as in \cite{RR} and
\cite{JV}, then assemble them to obtain a global biholomorphism from $M$ to $\cn$.

We begin by fixing some notation. By Proposition 3.1, let $0\le\lambda_1(t)\le
\lambda_2(t)\le...\le\lambda_n(t)$ be the eigenvalues of $R_\ijb(t)$ with respect
to $g(t)$ and let $0\le\mu_1<\mu_2\dots<\mu_l$ be their limits. Let $\rho>0$ and
$E_k(t)$, $1\le k\le l$ be as in Proposition 3.1 and let $P_k(t)$ be the orthogonal
projection onto $E_k(t)$ with respect to $g(t)$. Let $\tau_k=e^{-(\mu_k+1)T}$,
$1\le k\le l$. Note that for convenience, we have reversed the order of
$\lambda_i$ and hence the order of $\mu_k$.

By Theorem \ref{ns4t1}, $T_p^{(1,0)}(M)$ can be decomposed orthogonally with
respect to the initial metric as $E_1\oplus \dots\oplus E_l$ such that if  $v\in
E_k$ and $w\in E_j$ are nonzero vectors and if  $v(t)=v/|v|_t$, $w(t)=w/|w|_t$
where $|\cdot|_t$ is the norm taken with respect to $g(t)$,
 then
\begin{equation}\label{s5e4}
\lim_{t\to\infty}|P_k(t)v(t)|_t=1, \text{\ for  $1\le k\le l$ and\ }
\lim_{t\to\infty}\langle v(t),w(t)\rangle_t=0\text{\ for all $j\neq k$}.
\end{equation}
where $\langle\cdot,\cdot\rangle_t$ is the inner produce with respect to $g(t)$.
Moreover, the convergence is uniform.

For any $i$, let $E_{i,k}= d\Phi_{iT}^{-1}(E_k)$, $1\le k\le l$. Denote
$A(i)=A_i\cdots A_1$ and $A(i+j,i)=A_{i+j}\cdots A_{i+1}$. Then
$E_{i,k}=A(i)(E_{1,k})$ and $A_{i+1}(E_{i,k})=E_{i+1,k}$.
\begin{lem}\label{s5l1}
 Given $\e>0$, there exists $i_0$ such that if
$i\ge i_0$, then the following are true:
\begin{enumerate} \item[(i)]
$(1-\e)\tau_k\,|v|^2\le |A_{i}(v)|^2\le (1+\e)\tau_k\,|v|^2$ for all $v\in
E_{i,k}$ and $1\le k\le l$, where $\tau_k=e^{-(\mu_k+1)T}$. \item[(ii)] For any nonzero vector $v\in \C^n$
$$
(1-\e)\le\frac{|v|^2}{\sum_{k=1}^l|v_k|^2}\le (1+\e)
$$
where $v=\sum_{k=1}^l v_k$ is the decomposition of $v$ in
$E_{i,1}\oplus\cdots\oplus E_{i,l}$.
\end{enumerate}
\end{lem}
\begin{proof} (i) Let $1\le k\le l$. By (5.4), given $\e>0$, there exists $t_0>0$
such that
$$ |P_k(t)(w)|_t\geq 1-\e$$  for all $w\in E_k\setminus\{0\}$ and for all $t\ge t_0$. By the
definition of $E_k$ and  Proposition \ref{s3p1}, we have
that $|\Ric(w(t),\bar w(t))-\mu_k |\leq \e$ for all $w\in E_k\setminus \{0\}$,
provided $t_0$ is large enough. Suppose
$i_0>t_0/T$. Then for $i\geq i_0$ and $v\in E_{i,k}\setminus \{0\}$, there is
$w\in E_k\setminus\{0\}$ with $d\Phi^{-1}_{iT}(w)=v$. Hence
$A_{i}(v)=d\Phi^{-1}_{(i+1)T}(w)$. By (\ref{s5e1}), $|v|=|w|_{iT}$ and
$|A_{i}(v)|=|w|_{(i+1)T}$. By the K\"ahler-Ricci flow equation we have
\begin{equation}\nonumber
\begin{split}
\log\lf[\frac{|A_{i}(v)|^2}{|v|^2}\ri]+(\mu_k+1)T&= \log\lf[\frac{|w|^2_{(i+1)T}}{|w|^2_{iT}}\ri]+(\mu_k+1)T \\
&=\int_{iT}^{(i+1)T}\lf(\mu_k-\Ric(w(t),\bar w(t))\ri) dt.
\end{split}
\end{equation}
Since $|\Ric(w(t),\bar{w}(t))-\mu_k|\leq \e$ and $T$ is fixed, it is easy to see
that (i) is true.

(ii) Let $v\in \C^n$ be nonzero and let $v=\sum_{k=1}^lv_k$ be the
decomposition of $v$ in $E_{i,1}\oplus\cdots\oplus
E_{i,l}$. Let $w\in T_p^{(1,0)}(M)$ be such that
$d\Phi^{-1}_{iT}(w)=v$ and similarly decompose $w=\sum_{k=1}^l
w_k$ with respect to $E_1\oplus\cdots\oplus E_l$. Then
$v_k=d\Phi^{-1}_{iT}(w)$. Since $\langle v_j,v_k\rangle=\langle
w_j,w_k\rangle_{g(iT)}$ and $|v|^2=|w|^2_{g(iT)}$ by (\ref{s5e1}), (ii) follows from
(5.4).
\end{proof}

Let us fix more notation. Let $\Phi$ be a polynomial maps from $\C^n$ into
$\C^n$, which means that each component of $\Phi$ is a polynomial. Suppose $\Phi$ is of homogeneous of degree $m$. That is to say, each
component of $\Phi$ is a homogeneous polynomial of degree $m\ge 1$. We define
$$
||\Phi||=\sup_{v\in \C^n,v\neq0}\frac{|\Phi(v)|}{|v|^m}.
$$
In general, if $\Phi$ is a polynomial map with $\Phi(0)=0$, let $\Phi=\sum_{m=1}^q
\Phi_m$ be the decomposition of $\Phi$ such that $\Phi_m$ is homogeneous of
degree $m$, then $||\Phi||$ is defined as
$$
||\Phi||=\sum_{m=1}^q||\Phi_m||.
$$
If we decompose $\C^n$ as $E_{i,1}\oplus\cdots\oplus E_{i,l}$, we will denote
$\C^n$ by $\C^n_i$. Let $\Phi:\C_i^n\to \C_{i+1}^n$ be a map. Then we decompose
$\Phi$ as
$\Phi(v)=\sum_{k=1}^l\Phi_k(v)=\Phi_1\oplus\cdots\oplus\Phi_l$ where $\Phi_k(v)\in E_{i+1,k}$.  Let
$\a=(\a_1,\dots,\a_l)$ be a multi-index such that $|\a|=\sum_{k=1}^l\a_k=m\ge1$.
Then a polynomial map $\Phi$ is said to be {\it homogeneous of degree $\a$} if
$$
\Phi(c_1v_1\oplus\cdots\oplus c_lv_l)=c^\a\Phi(v_1\oplus\cdots\oplus v_l),
$$
where $v_k\in E_{i,k}$. Note that if  $\Phi$ homogeneous of degree $\a$, then $\Phi$ is homogeneous of degree $|\a|$ in the usual sense.  $\Phi$ is said to be {\it lower triangular}, if
$\Phi_k(v_1\oplus \cdots\oplus v_l)=c_kv_k+\Psi_k(v_1\oplus\cdots\oplus
v_{k-1})$.
\begin{lem}\label{s5l2}
 Let $\Phi:\C_i^n\to \C_{i+1}^{n}$ be
homogeneous of degree $\a=(\a_1,\dots,\a_l)$ with $|\a|=m$. Then
$$
 |\Phi(v_1\oplus\cdots\oplus v_l)| \le l^m||\Phi|| \,|v_1|^{\a_1}\cdots|v_l|^{\a_l}.
$$
Here by convention if $\a_i=0$, then $|v_i|^{\a_i}=1$ for all $v_i$.
 \end{lem}
\begin{proof} Let $v=v_1\oplus\cdots\oplus v_l$ such that $|v_k|=1$ for all $1\le
k\le l$, then
$$
|\Phi(v)|\le ||\Phi|| \, |v|^m\le l^m||\Phi||.
$$
Hence if $v_k\neq0$ for all $k$, then
$$
|\Phi(v)|=|\Phi(|v_1|\frac{v_1}{|v_1|}\oplus\cdots\oplus |v_l|\frac{ v_l}{
|v_l|})|\le l^m||\Phi|| \,|v_1|^{\a_1}\cdots|v_l|^{\a_l}.
$$
From this the lemma follows.
\end{proof}

Note that $\tau_1>\cdots>\tau_l$. Choose $1>\e>0$ small enough such that
$b^2(1-\e)^{-1}(1+\e)<1$ where $b<1$ is the constant in (5.3). Since we are
interested in the maps $F_i$ for large $i$, without loss of generality, we assume
the conclusions of Lemma 5.1 are true for all $i$ with this $\e$.
Let $m_0\ge2$ be a positive integer such that $a^{-1}b^{m_0}<\frac12$, where
$0<a<b<1$ are the constants in (\ref{s5e3}).

We now begin to assemble the maps $F_i$ to produce a
global biholomorphism from $M$ to $\C^n$.  The constructions follow those
in \cite{RR} and \cite{JV}; in particular those in \cite{JV} where the authors study
the dynamics of a randomly iterated sequence of biholomorphisms.

\begin{lem}\label{s5l3} Let $\Phi_{i+1}:\C_i^n\to \C_{i+1}^n$, $1\le i<\infty$, be a
family  homogeneous polynomial maps of degree $m\ge 2$ such that $\sup_i||
\Phi_i||<\infty$. Then  there exist  homogeneous polynomial maps $H_{i+1}$ and
$Q_{i+1}$ from $\C_i^n$ to $\C_{i+1}^n$    such that
$\Phi_{i+1}=Q_{i+1}+H_{i+1}-A_{i+2}^{-1}H_{i+2} A_{i+1}$. Moreover, $H_{i+1}$
and $Q_{i+1}$ satisfy the following: \begin{enumerate} \item[(i)]
$\sup_{i}||H_i||<\infty$ and $\sup_i||Q_i||<\infty$. \item[(ii)]  $Q_{i+1}=0$ if
$m\ge m_0$. \item[(iii)] $Q_{i+1}$ is lower triangular:
$$Q_{i+1}(v_1\oplus \cdots\oplus v_l)=0\oplus Q_{i+1,2}(v_1)\oplus Q_{i+1,3}
(v_1\oplus v_2)\oplus\cdots \oplus Q_{i+1,l}(v_1\oplus\cdots\oplus v_{l-1})$$
where $v_k\in E_{i,k}$ and $Q_{i+1,k}:\C_i^n\to E_{i+1,k}$.\end{enumerate}
\end{lem}
\begin{proof} For each $i$, let $\beta_k$ be a unitary basis for $E_{i,k}$ with
respect to the standard metric of $\C_i^n$. Let $v\in \C_i^n$ and if
$v=\sum_{k=1}^l\sum_{w\in\beta_k}a_w w$, then
$$C_1^{-1}|v|^2\le\sum_{k=1}^l \sum_{w\in \beta_k}|a_w|^2\le C_1|v|^2
$$
for some constant $C_1$ independent of $i$ by Lemma 5.1(ii). Hence if we
decompose $\Phi_{i+1}$ into $\a$-homogeneous parts $\Phi_{i+1,\a}$, $|\a|=m$,
then $||\Phi_{i+1,\a}||\le C_2||\Phi_{i+1}||$ for some constant $C_2$ independent
of $\Phi_{i+1}$ and $i$. Moreover, if we decompose
$\Phi_{i+1,\a}=\Phi_{i+1,\a,1}\oplus\cdots\oplus \Phi_{i+1,\a,l}$ with
$\Phi_{i+1,\a,k}(v)\in E_{i+1,k}$, then by Lemma 5.1(ii) again,
$||\Phi_{i+1,\a,k}||\le C_3||\Phi_{i+1,\a}||$ for some constant $C_3$ independent
of $i$. Hence in order to prove the lemma, we may assume that $\Phi_{i+1}$ is
homogeneous of degree $\a=(\a_1,\dots,\a_l)$ with $|\a|=m$ and $\Phi_{i+1}(v)\in
E_{i+1,k}$ for all $i$ for some   $1\le k\le l$.

Suppose $m\ge m_0$.  Then we define $Q_{i+1}=0$ and let
$$
H_{i+1}=\Phi_{i+1}+\sum_{s=0}^\infty A_{i+2}^{-1}\cdots
A_{i+s+2}^{-1}\Phi_{i+s+2}A_{i+s+1}\cdots A_{i+1}.
$$
To see $H_{i+1}$ is well-defined, by (\ref{s5e3}) we have that for any $v\in \C_i^n$,
$$
|\Phi_{i+s+2}A_{i+s+1}\cdots A_{i+1}(v)|\le ||\Phi_{i+s+2}||(b^{s+1}|v|)^m
$$
and
\begin{equation}\nonumber
\begin{split}
|A_{i+2}^{-1}\cdots A_{i+s+2}^{-1}\Phi_{i+s+2}A_{i+s+1}\cdots A_{i+1}(v)|&\le ||\Phi_{i+s+2}||\lf(a^{-1}b^m\ri)^{s+1} |v|^m\\
&\le 2^{-s-1}||\Phi_{i+s+2}||\,|v|^m.
\end{split}
\end{equation}
Hence $H_{i+1}$ is well-defined, homogeneous of degree $m$ and $||H_{i+1}||\le
C_4$ for some constant $C_4$ independent of $i$. It is easy to see  that
$H_{i+1}$ and $Q_{i+1}$ satisfy the required conditions.

Now consider the case that $2\le m<m_0$. Decompose $\Phi_{i+1}$ as
$\Phi_{i+1}^{(1)}+\Phi_{i+1}^{(2)}$ where $\Phi_{i+1}^{(1)}(v_1\oplus\cdots
\oplus v_l)=\Phi_{i+1}(v_1\oplus\cdots\oplus v_{k-1}\oplus 0\cdots\oplus 0)$
consisting all terms that depending only  on $v_1,\dots, v_{k-1}$ and
$\Phi_{i+1}^{(2)}=\Phi_{i+1}-\Phi_{i+1}^{(1)}$. Let $Q_{i+1}=\Phi_{i+1}^{(1)}$.
Since $\Phi_{i+1}(v)\in E_{i+1,k}$, it is easy to see that $Q_{i+1}$ satisfies
condition (iii) in the lemma. It is also easy to see that $||Q_{i+1}||\le
||\Phi_{i+1}||$.

Suppose  $\a_j=0$ for all $j\ge k$, then $\Phi^{(2)}_i=0$. In this case, let
$H_{i+1}=0$. Then $Q_{i+1}$ and $H_{i+1}$ satisfy the required conditions.

Suppose there is $j\ge k$ with $\a_j\ge 1$.  Then define
\begin{equation}\label{s5e5}
H_{i+1}=\Phi_{i+1}^{(2)}+\sum_{s=0}^\infty A_{i+2}^{-1}\cdots
A_{i+s+2}^{-1}\Phi^{(2)}_{i+s+2}A_{i+s+1}\cdots A_{i+1}.
\end{equation}
To prove $H_{i+1}$ is well-defined and $||H_{i+1}||$ is uniformly bounded, we
observe that
\begin{equation}\label{s5e6}
||\Phi^{(2)}_{i+s+2}||\le ||\Phi_{i+s+2}^{(1)}|| +||\Phi_{i+s+2}||\le
2||\Phi_{i+s+2}||.
\end{equation}
Let $v\in \C^n_i$ and let $w=w_1\oplus \cdots\oplus w_l=A(i+s+1,i)(v)$ and let
$u=A(i+s+2,i+1)^{-1}(\Phi_{i+s+2}^{(2)}(w))$. Note that if $v=v_1\oplus\cdots\oplus
v_l$ with $v_q\in E_{i,q}$, then $A_{i+r}(v_q)\in E_{i+r,q}$. Hence  by Lemma
\ref{s5l2},   Lemma \ref{s5l1}(i) and (\ref{s5e3})
\begin{equation}\nonumber
\begin{split}
|\Phi_{i+s+2}^{(2)}(w)|&\le l^m ||\Phi^{(2)}_{i+s+2}|| \,|w_1|^{\a_1}\dots|w_l|^{\a_l} \\
&\le 2l^m||\Phi_{i+s+2}|| |w|^{m-1}\,|w_j|\\
& \le 2l^m||\Phi_{i+s+2}||  b^{(s+1)(m-1)}\lf[(1+\e)    \tau_j\ri]^{\frac
{s+1}2}\,|v|^m.
\end{split}
\end{equation}
Since $\Phi_{i+s+2}^{(2)}w\in E_{i+s+2, k}$, by Lemma \ref{s5l1}(i) and the fact that
$A_{r+1}^{-1}(E_{r+1,k})=E_{r,k}$ for all $r$, we have
\begin{equation}\label{s5e7}
\begin{split}
|u|&=|A(i+s+2,i+1)^{-1}(\Phi_{i+s+2}^{(2)}w)|\\
& \le\lf[(1-\e)\tau_k\ri]^{-\frac{s+1}2}|\Phi_{i+s+2}^{(2)}w |\\
&\le 2l^m||\Phi_{i+s+2}||  b^{(s+1)(m-1)}\lf[(1-\e)\tau_k\ri]^{-\frac{s+1}2}\,\lf[(1+\e)    \tau_j\ri]^{\frac {s+1}2}\,|v|^m\\
&\le 2l^m||\Phi_{i+s+2}||\lf[b^2(1-\e)^{-1}(1+\e)\ri]^{\frac{s+1}2}
\end{split}
\end{equation}
since $\tau_k\ge \tau_j$ for $j\ge k$, $m\ge 2$ and $b<1$. Since we have chosen
$\e$ such that $b^2(1-\e)^{-1}(1+\e)<1$, from (\ref{s5e5})--(\ref{s5e7}), we conclude that
$H_{i+1}$ is well-defined and $||H_{i+1}||$ are uniformly bounded. Note that
$H_{i+1}$ is homogeneous of degree $m$. Then $Q_{i+1}$ and $H_{i+1}$ satisfy the
required conditions.
\end{proof}
\begin{lem}\label{s5l4} Given any $m\ge2$, we can find constants $C(m)>0$ and
$r_1\ge r_m>0$ and families of holomorphic maps $T_{i,m}$ from $D(r_m)\subset
\C_{i}^n$ to $D(r_m)\subset \C_{i}^n$ and $G_{i+1,m}$ from $ \C_i^n$ to
$\C_{i+1}^n$ with the following properties: \begin{enumerate} \item[(i)] For each
$i$, $T_{i+1,m}$ is a polynomial map of degree $m-1$ which is biholomorphic to
its image, $T_{i+1,m}(0)=0$, $T'_{i+1,m}(0)=Id$ and $||T_{i+1,m}||\le C(m)$.
\item[(ii)] $G_{i+1,m}=A_{i+1}+\wt G_{i+1,m}$ where $\wt G_{i+1,m}$ is a
polynomial map of degree $m-1$,
$$\wt G_{i+1,m}(v_1\oplus\cdots\oplus v_l)=0\oplus \wt G_{i+1,m, 2}(v_1)\oplus\cdots\oplus \wt G_{i+1,m, 2}(v_1\oplus\cdots\oplus v_{l-1})$$ is lower triangular, and $||G_{i+1,m}||\le C(m)$, $\wt G_{i+1,m}(0)=0$ and $G'_{i+1m}(0)=0$. Moreover, $G_{i+1,m}=G_{i+1,m_0}$ for all $m\ge m_0$, where $m_0$ is the integer in Lemma 5.3.
\item[(iii)] $F_{i+1}(D(r_m))\subset D(r_m)$ and
$$
|T_{i+1,m}F_{i+1}(v)-G_{i+1,m}T_{i,m}(v)|\le C(m)|v|^m.
$$
\end{enumerate}
Here $T_{i+1,m}F_{i+1}-G_{i+1,m}T_{i,m}$ means $T_{i+1,m}\circ F_{i+1}
-G_{i+1,m}\circ T_{i,m}$.
\end{lem}

\begin{proof} Note that since $A_{i+1}$ is nonsingular, $G_{i+1,m}$ will be a
biholomorphism. We will construct the maps by induction. For $m=2$, let
$T_{i+1,m}=Id$, $G_{i+1,m}=A_{i+1}$. Since $F_{i+1}(D(r_1))\subset D(r_1)$ and is
holomorphic, by (\ref{s5e2}) we can take $r_2=\frac12 r_1$, then it is easy to see that one can
find $C_2$ satisfies the required conditions. Suppose we have found $T_{i+1,m}$,
$G_{i+1,m}$, $C(m)$ and $r_m$ which have the required properties. Since
$$
|T_{i+1,m}F_{i+1}(v)-G_{i+1,m}T_{i,m}(v)|\le C(m)|v|^m
$$
we have $||\Phi_{i+1}||\le C_1$ for some $C_1$ which is independent of $i$, where
$\Phi_{i+1}$ is the homogeneous polynomial of degree $m$ which is the $m$-th
power terms of the Taylor series of $T_{i+2,m}F_{i+1}-G_{i+1,m+1}T_{i+1,m}$. By
Lemma 5.3, we can find $H_{i+1}$ and $Q_{i+1}$ such that both are homogeneous of degree
$m$, $H_{i+1}$ and $Q_{i+1}$ satisfies conditions (i)--(iii) in Lemma
\ref{s5l3} and
$$
\Phi_{i+1}=Q_{i+1}+H_{i+1}-A_{i+2}^{-1}H_{i+2}A_{i+1}.
$$
Now define $T_{i,m+1}=T_{i,m}+A_{i+1}^{-1}H_{i+1}$ and
$G_{i+1,m+1}=G_{i+1,m}+Q_{i+1}$. Note that if $m\ge m_0$, then $Q_{i+1}=0$. By
the induction hypothesis, Lemma \ref{s5l3} and (\ref{s5e3}), it is easy to see that
$T_{i+1,m+1}$ and $G_{i+1,m+1}$ satisfy (i) and (ii) of the lemma for some
constants $C(m+1)$ and $r_{m+1}\le \frac12 r_m$. It remains to check condition
(iii). We proceed as in \cite{RR}.

 In the following, $O(m+1)$ will denote some function $h$
 such that $|h(v)|\le C|v|^{m+1}$ for $|v|\le \frac12 r_m$,
 where $C$ is a constant independent of $i$.
\begin{equation}\label{s5e8}
\begin{split}
&T_{i+1,m+1}F_{i+1}-G_{i+1,m+1}T_{i,m+1}\\
&\quad =(T_{i+1,m}+A_{i+2}^{-1}H_{i+2})F_{i+1}-(G_{i+1,m}+Q_{i+1})(T_{i,m}+A_{i+1}^{-1}H_{i+1})\\
&\quad =\lf[T_{i+1,m}F_{i+1}-G_{i+1,m}T_{i,m}\ri]+G_{i+1,m}T_{i,m}-G_{i+1,m}(T_{i,m}+A_{i+1}^{-1}H_{i+1})\\
&\qquad -Q_{i+1} (T_{i,m}+A_{i+1}^{-1}H_{i+1})+A_{i+2}^{-1}H_{i+2}F_{i+1}
\end{split}
\end{equation}
Since $F_i(D(r_m))\subset D(r_m)$, and $||T_{i,m}||$ and $||G_{i,m}||$ are uniformly
bounded,
\begin{equation}\nonumber
\begin{split}
T_{i+1,m}F_{i+1}-G_{i+1,m}T_{i,m}&=\Phi_{i+1}+O(m+1)\\
&=Q_{i+1}+H_{i+1}-A_{i+2}^{-1}H_{i+2}A_{i+1}+O(m+1).
\end{split}
\end{equation}
Combining this with (5.8), we have
\begin{equation}\label{s5e9}
\begin{split}
&T_{i+1,m+1} F_{i+1}-G_{i+1,m+1}T_{i,m+1}\\
&=Q_{i+1}
+H_{i+1}-A_{i+2}^{-1}H_{i+2}A_{i+1}+G_{i+1,m}T_{i,m}-G_{i+1,m}(T_{i,m}+A_{i+1}^{-1}H_{i+1})\\
&\quad   -Q_{i+1} (T_{i,m}+A_{i+1}^{-1}H_{i+1})+A_{i+2}^{-1}H_{i+2}F_{i+1}+O(m+1)\\
&=\lf[G_{i+1,m}T_{i,m}-G_{i+1,m}(T_{i,m}+A_{i+1}^{-1}H_{i+1})+H_{i+1}\ri]\\
&\quad +\lf[Q_{i+1}-Q_{i+1}
(T_{i,m}+A_{i+1}^{-1}H_{i+1})\ri]+\lf[A_{i+2}^{-1}H_{i+2}F_{i+1}-A_{i+2}^{-1}H_{i+2}A_{i+1}\ri]\\
&\quad +O(m+1).
\end{split}
\end{equation}
Denote the differential of a map $h$ by $h'$. Then
\begin{equation}\nonumber
\begin{split}
H_{i+2}\circ F_{i+1}-H_{i+2}A_{i+1}&=\int_0^1\frac{d}{ds}(H_{i+1}\circ (sF_{i+1}-(1-s)A_{i+1})ds\\
&=\int_0^1\lf[ H_{i+1}'(sF_{i+1}-(1-s)A_{i+1})\ri](F_{i+1}-A_{i+1})ds
\end{split}
\end{equation}
where the multiplication of the terms under integral sign is matrix
multiplication. By (\ref{s5e2}), (\ref{s5e3}), the definition of $A_{i+1}$ and the fact that
$||H_{i+1}||$ are uniformly bounded and homogeneous of degree $m\ge 2$, we have
\begin{equation}\label{s5e10}
H_{i+2}\circ F_{i+1}-H_{i+2}A_{i+1}=O(m+1).
\end{equation}
Using (5.3) and the facts that $||Q_{i+1}||$, $||T_{i,m}||$ and $||H_{i+1}||$ are
uniformly bounded, $Q_{i+1}$ is homogeneous of degree $m\ge 2$ and that $T'_{i,m}(0)=Id$, we can prove similarly that
\begin{equation}\label{s5e11}
Q_{i+1}-Q_{i+1} \circ(T_{i,m}+A_{i+1}^{-1}H_{i+1})=O(m+1).
\end{equation}
 Finally,
\begin{equation}\nonumber
\begin{split}
&G_{i+1,m}\circ T_{i,m}-G_{i+1,m}\circ(T_{i,m}+A_{i+1}^{-1}H_{i+1})+H_{i+1}\\
&\quad  =-\int_0^1\frac{d}{ds}\lf(G_{i+1,m}\circ(T_{i,m}+sA_{i+1}^{-1}H_{i+1})\ri)ds+H_{i+1}\\
&\quad =-\int_0^1 \lf(\lf[G_{i+1,m}'(T_{i,m}+sA_{i+1}^{-1}H_{i+1})\ri](A_{i+1}^{-1}H_{i+1})-A_{i+1}A_{i+1}^{-1}H_{i+1}\ri) ds\\
&\quad =-\int_0^1
\lf(\lf[G_{i+1,m}'(T_{i,m}+sA_{i+1}^{-1}H_{i+1})-A_{i+1}\ri](A_{i+1}^{-1}H_{i+1})\ri)
ds
\end{split}
\end{equation}
Using (\ref{s5e3}) and the facts that $G_{i+1,m}'(0)=A_{i+1}$, that $||G_{i+1,m}||$,
$||H_{i+1}||$ are uniformly bounded, and that $H_{i+1}$ is homogeneous of degree
$m$ we conclude that
\begin{equation}\label{s5e12}
G_{i+1,m}\circ
T_{i,m}-G_{i+1,m}\circ(T_{i,m}+A_{i+1}^{-1}H_{i+1})+H_{i+1}=O(m+1).
\end{equation}
From (\ref{s5e9})--(\ref{s5e12}), we conclude that
$$
|T_{i+1,m+1}F_{i+1}(v)-G_{i+1,m+1}T_{i,m+1}(v)|\le C(m+1)|v|^{m+1}.
$$
This completes the proof of the lemma.
\end{proof}

Let $m\ge m_0$ and denote $G_{i+1,m}$ simply by $G_{i+1}$ and denote $\wt
G_{i+1,m}$ by $\wt G_{i+1}$ etc. Note that $G_{i+1}$ is independent of $m$ and is
a biholomorphism on $\C^n$. Degree of each $G_{i+1}$ is $m-1$. For any positive
integers $i,j$, let $G(i+j,i)=G_{i+j}\cdots G_{i+1}$.
\begin{lem}\label{s5l5} Let
$G_{i+1}$ as above, then its inverse is a polynomial map of degree $(m-1)^{l-1}$
and satisfies:
$$
G_{i+1}^{-1}=A_{i+1}^{-1}+S_{i+1}
$$
where $S_{i+1}:\C^n_{i+1}\to \C^n_i$ with
$$S_{i+1}(w_1\oplus \cdots\oplus w_l)=0\oplus S_{i+1,2}(w_1)\oplus\cdots\oplus S_{i+1,l}(w_1\oplus \cdots\oplus w_{l-1}).$$
Moreover, $||G_{i+1}^{-1}||$ is bounded by a constant independent of $i$.
\end{lem}
\begin{proof} Let $w_1\oplus \cdots\oplus w_l=E_{i+1,1}\oplus \cdots\oplus
E_{i+1,l}=\C^n_{i+1}$. Let $v_1=A_{i+1}^{-1}w_1$, $v_2=A_{i+1}^{-1}(w_2-\wt
G_{i+1,2}(v_1)),\dots, v_l=A_{i+1}^{-1}(w_l-\wt G_{i+1,l}(v_1\oplus \cdots\oplus
v_{l-1}))$. Let $S_{i+1,k}(w_1\oplus \cdots\oplus w_{k-1})=-A_{i+1}^{-1}\wt
G_{i+1,k}(v_1\oplus \cdots\oplus v_{k-1})$, $2\le k\le l$. It is easy to see that
$S_{i+1,k}$ is well-defined and $S_{i+1,k}(w_1\oplus \cdots\oplus w_{k-1})\in
E_{i,k}$ because $A_{i+1}(E_{i,k})=E_{i+1,k}$. Moreover, the degree of each
$S_{i+1,k}$ is at most $(m-1)^{k-1}$. It is also easy to see that

$$
G_{i+1}^{-1}=A_{i+1}^{-1}+S_{i+1}
$$
where $S_{i+1}=0\oplus S_{i+1,2}\oplus\cdots\oplus S_{i+1,l}$.

Let $w_1\oplus \cdots\oplus w_l\in \C^n_{i+1}$ with $|w_k|\le 1$ and  $v_1\oplus
\cdots\oplus v_l=G_{i+1}^{-1}(w_1\oplus \cdots\oplus w_l)$. We claim that $|v_k|$
is bounded by a constant independent of $i$ for each $k$. If this is true, then
by Lemma \ref{s5l1} and (\ref{s5e3}) again, we can conclude that $||G_{i+1}^{-1}||$ is bounded
by a constant independent of $i$. To prove the claim, by (\ref{s5e3}) have
$|v_1|=|A_{i+1}^{-1}(w_1)|$ is uniformly bounded for $|w_1|\le 1$. Since
$||G_{i+1}||$ is uniformly bounded by a constant independent of $i$, $||\wt
G_{i+1,k}||$ is also uniformly bounded by a constant independent of $i$ by Lemma
\ref{s5l1}(ii) and (\ref{s5e3}).  Suppose we have proved that $|v_1|,\dots,|v_{k-1}|$ are
bounded by a constant independent of $i$. Then
$$|S_{i+1,k}(w_1\oplus\cdots\oplus w_{k-1})| =|A_{i+1}^{-1} \wt G_{i+1,k}(v_1\oplus\cdots\oplus w_{k-1})|
$$
is also bounded by a constant independent of $i$. Hence $|v_k|$ is bounded by a
constant independent of $i$. This completes the proof of the lemma.
\end{proof}

\begin{lem}\label{s5l6} Let $D(R)$ be a ball in $\C^n$ with radius $R$ with center
at the origin.  Then the following are true: \begin{enumerate} \item[(i)] There exist
$\beta>0$ such that for all $z, z'\in D(R)$ and for any positive integers $i$ and
$j$,
$$
|G(i+j,i)^{-1}(z)-G(i+j,i)^{-1}(z')|\le \beta^j |z-z'|.
$$
\item[(ii)] For any positive integer $i$ and for any open set $U$ containing the
origin,
$$
\bigcup_{j=1}^\infty G(i+j,i)^{-1}(U)=\C^n.
$$
\end{enumerate}
\end{lem}
\begin{proof}  (i) For simplicity, let us assume that $R=1$ and  first assume that
$i=0$.   Let us write
\begin{equation}\label{s5e13}
G(j,0)^{-1}=G_1^{-1}\cdots G_j^{-1} =H_{j,1}\oplus\cdots \oplus H_{j,l}
\end{equation} with $H_{j,k}(v)\in E_{1,k}$. By Lemma 5.2 and the Schwartz
lemma, it is sufficient to prove that
$$|H_{j,k}(v)|\le \beta^j
$$ for some constant $\beta$ for all $k$ and $j$ provided $|v|\le 1$. By Lemma 5.5, $G_i^{-1}=A_i+S_i$ where $S_i$ satisfies the conclusions in the lemma. Let $v=v_1\oplus\cdots\oplus v_l\in \C_{j}^n$. Then $H_{j,1}(v)=A_{j}^{-1}\cdots A_1^{-1}(v_1)\le a^j|v_1|\le 2a^j$ , where we have used Lemma 5.2. Hence (5.13) is true for $k=1$. Suppose (\ref{s5e13}) is true for $1,\dots,k-1$. We may assume that $\beta>a^{-1}$. By Lemma \ref{s5l2} and \ref{s5l5}, we know that $||S_j||$ is uniformly bounded. Let $C_j=\max_k\{\max_{|v|\le 1}|H_{j,k}(v)|,1\}$. Since $G_{j,k}(w)=A_j^{-1}(w_k)+S_j(w_1\oplus \cdots \oplus w_{k-1})$, we have
\begin{equation}\nonumber
\begin{split}
C_j&\le a^{-1}C_{j-1}+C\beta^{(j-1)N}\\
& \le 2C_{j-1}\beta_1^{j-1}
\end{split}
\end{equation}
where $N=(m-1)^{l-1}$ which is the degree of $S_i$,   $C>\ge1$ is a constant
dependent only on $||S_i||$ and $N$, and $\beta_1=C\beta^N\ge a^{-1}$, where we
have used the fact that $C_{j-1}\ge 1$. Hence $C_j\le (2\beta_1)^{j-1}C_1$.  From
this the lemma follows for $i=0$. For general $i$, the proof is similar. Note
that the constants in the proof do not depend on $i$.

(ii) The proof is similar to the proof of (i). Let us write $G_j\cdots
G_1=K_{j,1}\oplus\cdots\oplus K_{j,l}$. Then $K_{j,1}(v_1\oplus\cdots\oplus
v_l)=A_j\cdots A_1(v_1)$. Hence $K_{j,1}(v)$ converge to zero uniformly on
compact sets. Suppose $K_{j,1},\dots,K_{j,k-1}$ converge uniformly to 0 on
compact sets. Let $\Omega$ be a compact set and let $s_j=\sup_{v\in
\Omega}|K_{j,k}$. Then as before,
$$
s_j\le b s_{j-1}+\sup_{v\in \Omega}|\wt
G_{j,k}(K_{j-1,1}(v),\dots,K_{j-1,k-1}(v))|.
$$
Hence
$$\limsup_{j\to\infty}s_j\le b\limsup_{j\to\infty} s_{j-1}
$$
because $||\wt G_{j,k}||$ are uniformly bounded with uniformly bounded degrees
and $K_{j-1,p}(v)\to0$ uniformly on $\Omega$ for $1\le p\le k-1$. From this it is
easy to see that $s_j\to0$ as $j\to\infty$. Hence $G_j\cdots G_1\to0$ uniformly
on compact sets. From this (ii) follows.
\end{proof}

Let $\beta$ be the constant in Lemma \ref{s5l6}. Note that $\beta$ does not depend on
$i$ and $m$ provided $m\ge m_0$, where $m_0\ge 2$ is the integer in Lemma 5.3.
Fix $m\ge m_0$ such that
\begin{equation}\label{s5e14}
\delta^m\le \frac12 \beta.
\end{equation}
where $1>\delta>0$ be the constant in (\ref{s5e2}). Let $G_{i,m}$, $T_{i,m}$   be the
maps given in Lemma 5.4 which are defined on $D(r_m)$, $0<r_m<r_1<1$. Let us
denote $G_{i,m}$ by $G_i$, $T_{i,m}$ by $T_i$ and $r_m$ be $r$.

In the following, a holomorphic map $\Phi$ from a complex manifold to another is
said to be {\it nondegenerate} if it is injective and so that it is a
biholomorphism onto its image. We apply the method in \cite{RR} to obtain the following.

\begin{lem}\label{s5l7} Let $k\ge0$ be an integer. Then
$$
\Psi_k=\lim_{l\to\infty}G^{-1}_{k+1}\circ G^{-1}_{k+2}\circ\cdots\circ
G^{-1}_{k+l}\circ T_{k+l}\circ F_{k+l}\circ\cdots\circ F_{k+2}\circ F_{k+1}
$$
exists and is a nondegenerate holomorphic map from $D(r)$ to $\C^n$. Moreover,
there is a constant $\gamma>0$ which is independent of $k$ such that
\begin{equation}\label{s5e15}
\gamma^{-1}D(r)\subset\Psi(D(r))\subset \gamma D(r).
\end{equation}
\end{lem}
\begin{proof} Let $\Phi_l=G^{-1}_{k+1}\circ G^{-1}_{k+2}\circ\cdots\circ
G^{-1}_{k+l}\circ T_{k+l}\circ F_{k+l}\circ\cdots\circ F_{k+2}\circ F_{k+1}. $ By
the construction in Lemma \ref{s5l4}, $\Phi_l$ is a nondegenerate holomorphic map on
$D(r)$ and $\Phi_l(0)=0$. For any $z\in D(r)$, let $w=F_{k+l}\circ\cdots\circ
F_{k+1}(z)$. Then $|w|\le \delta^l r$ by (\ref{s5e2}). Hence $T_{k+l}(w)\in D(R)$,
$T_{k+l+1}\circ F_{k+l+1}(w)\in D(R)$, $G^{-1}_{k+l+1}\circ T_{k+l+1}\circ
F_{k+l+1}(w)\in D(R)$, and $G_{k+l+1}\circ T_{k+l}(w)\in D(R)$ for some $R$
independent of $k$  and $l$ by Lemmas \ref{s5l4} and \ref{s5l5}. By Lemmas \ref{s5l4}(iii) and \ref{s5l6}, we
have
\begin{equation}\nonumber
\begin{split}
\bigg|G&^{-1}_{k+1} \circ\cdots\circ G^{-1}_{k+l}\circ  G^{-1}_{k+l+1}\circ T_{k+l+1}\circ F_{k+l+1}\circ F_{k+l}\cdots \circ F_{k+1}(z)\\
&\qquad -
G^{-1}_{k+1}\circ  \cdots\circ G^{-1}_{k+l}\circ T_{k+l}\circ F_{k+l}\circ\cdots \circ F_{k+1}(z)\bigg|\\
&=\bigg|G^{-1}_{k+1} \circ\cdots\circ G^{-1}_{k+l}\circ  G^{-1}_{k+l+1}\circ
T_{k+l+1}\circ F_{k+l+1}(w)\\
&\qquad -G^{-1}_{k+1}\circ  \cdots\circ G^{-1}_{k+l}\circ T_{k+l}(w)\bigg|\\
&\le \beta^l\bigg|G^{-1}_{k+l+1}\circ T_{k+l+1}\circ F_{k+l+1}(w)-T_{k+l}(w)\bigg|\\
&\le \beta^{l+1}\bigg|T_{k+l+1}\circ F_{k+l+1}(w)-G_{k+l+1}\circ T_{k+l}(w)\bigg|\\
&\le C_1\beta^{l+1}|w|^m\\
&\le C_1\beta^{l+1}\delta^{lm}\\
&\le C_1\beta\lf(\frac12\ri)^l
\end{split}
\end{equation}
by (\ref{s5e14}). From this it is easy to see that $\Psi_k=\lim_{l\to\infty}\Phi_l$
exists and is holomorphic on $D(r)$. Moreover
$$|\Psi_k(z)|\le |\Phi_1(z)|+C_1\beta.$$
Using   (\ref{s5e3}) and the fact that $||G_i||$ and $||T_i||$ are uniformly bounded,
we can find $\gamma>1$ independent of $k$ and $l$ such that $\Psi_k(D(r))\subset
\gamma D(r)$.  Since $\Phi_l'(0)=Id$, $\Psi_k'(0)=Id$. By the gradient estimates
of holomorphic functions, $|\Phi_k'(z)-Id|\le C_2|z|$ on $\frac12 D(r)$ for some
constant $C_2$ independent of $k$. Hence there exists $r>r'>0$ independent of $k$
such that $\Phi_k$ is nondegenerate in $D(r')$ and $\Psi_k(D(r))\supset
\gamma^{-1}D(r)$ provided $\gamma$ is large enough independent of $k$. To prove
that $\Psi_k$ is nondegenerate on $D(r)$, let $l_0$ be such that
$F_{k+l_0}\cdots\circ F_{k+1}(D(r))\subset D(r')$. Then
$$
\Psi_k=G^{-1}_{k+1}\circ\cdots\circ G^{-1}_{k+l_0}\circ \Psi_{k+l_0}\circ
F_{k+l_0}\circ\cdots F_{k+1}.
$$
Since $F_{k+l_0}\cdots\circ F_{k+1}$ is nondegenerate on $D(r)$, $\Psi_{k+l_0}$
is nondegenerate on $D(r')$, and $G^{-1}_{k+1}\circ\cdots\circ G^{-1}_{k+l_0}$ is
a biholomorphism of $\C^n$, we conclude that $\Psi_k$ is nondegenerate on $D(r)$.
\end{proof}
Now we are ready to prove the following uniformization theorem.
\begin{thm}\label{s5t1} Let $(M^n,\wt g)$ be a complete noncompact K\"ahler manifold
with nonnegative and bounded holomorphic bisectional curvature. Suppose the
scalar curvature of $M$ satisfies
\begin{equation}\label{s5e16}
\frac{1}{V_x(r)}\int_{B_x(r)}R\le \frac{C}{1+r^2}
\end{equation}
for some constant $C$ for all $x\in M$ for all $r$. Suppose $(M, g)$ has
maximal volume growth. Then $M$ is biholomorphic to $\C^n$.  Moreover, the
assumption of maximal volume growth can be removed if $M$ has positive curvature
operator.
\end{thm}
\begin{proof} If $\wt g$ satisfies the given conditions, then one can
  solve the
K\"ahler-Ricci flow (\ref{krfn}) and construct $\Phi_t$ and $F_i$ as in the beginning of this
section. We can also construct $G_i$, $T_i$ as in Lemma \ref{s5l3} so that Lemmas \ref{s5l6}
and \ref{s5l7} are true. Let $\Omega_i=\Phi_{iT}^{-1}(D(r))$ where $r>0$ is the constant
in Lemma \ref{s5l7}. By (\ref{s5e1}) and the fact that the solution $g(t)$ of (\ref{krfn})
decays exponentially, $\{\Omega_i\}_{i\ge 1}$ exhausts $M$. Consider the following
holomorphic maps from $\Omega_{i}$ to $\C^n$:
$$
S_i=G^{-1}_1\circ\cdots G^{-1}_{i}\circ T_{i}\circ \Phi_{iT}^{-1}.
$$
For each fixed $k$, and $l\ge 1$
\begin{equation}\nonumber
\begin{split}
S_{k+l} &=G^{-1}_1\circ\cdots \circ G^{-1}_k\circ G^{-1}_{k+l}\circ T_{k+l}\circ \Phi_{(k+l)T}^{-1}\\
&=G^{-1}_1\circ\cdots \circ G^{-1}_k\circ\lf[ G^{-1}_{k+1}\circ\cdots\circ
G^{-1}_{k+l}\circ T_{k+l}\circ F_{k+l}\circ \cdots\circ F_{k+1}\ri]\circ
\Phi^{-1}_{kT}
\end{split}
\end{equation}
By Lemma \ref{s5l7}, we conclude that $S=\lim_{i\to\infty}S_i$ exists and is a
nondegenerate holomorphic map from $M$ into $\C^n$. Moreover,
$S=G^{-1}_1\circ\cdots \circ G^{-1}_k\circ \Psi_k\circ \Phi^{-1}_{kT}$ on
$\Omega_k$ where $\Psi_k$ is the nondegenerate holomorphic map in Lemma \ref{s5l7}.
Hence
$$
S(\Omega_k)=G^{-1}_1\circ\cdots \circ G^{-1}_k\circ \Psi_k(D(r))\supset G^{-1}_1\circ\cdots \circ G^{-1}_k(\gamma^{-1}D(r))
$$
by Lemma \ref{s5l7}, for some $\gamma$ independent of $k$. Therefore $S(M)=\C^n$ by
Lemma \ref{s5l6}(ii). This completes proof of the theorem.
\end{proof}
By a recent result of Ni \cite{Ni2}, if $M$ has maximal
volume growth, then (\ref{s5e16}) is satisfied automatically. Hence we have:
\begin{cor}\label{s5c1} Let $(M^n,\wt g)$ be a complete noncompact K\"ahler manifold
with nonnegative and bounded holomorphic bisectional curvature. Suppose $M$ has
maximal volume growth, then $M$ is biholomorphic to $\C^n$.
\end{cor}

We also have  the following
uniformization theorem. 
\begin{thm}\label{s5t2} Let $(M^n,\wt g)$ be a complete
noncompact K\"ahler manifold with nonnegative curvature operator such that the
scalar curvature $R$ of $M$ satisfies (5.16). Then the universal cover of $M$ is
biholomorphic to $\C^n.$
\end{thm}
\begin{proof}
  Let $\tilde{g}(t)$ be the corresponding solution to the K\"ahler-Ricci flow \ref{krf}. Let $\widetilde M$ be the universal cover of $M$. We then lift the flow $\tilde{g}(t)$
  to $\widetilde M$ and denote the lifted flow by $\tilde{h}(t)$. 

By the result in \cite{cao1}  and the De Rham decomposition theroem, 
one may assume that $\wt M=\C^{k}\times N_1\times \cdots \times N_l$
  isometrically and holomorphically so that each  $N_j$ 
is irreducible and  has nonnegative curvature operator and positive
  Ricci curvature. Note that the flow $\tilde{h}(t)$ still satisfies 
the K\"ahler-Ricci flow equation when restricted on each $N_j$.
Now suppose  there is a positive constant $C$ 
such that for $t$ large enough, the injectivity radius of $\tilde{h}(t)$ is bounded
below by $Ct^{1/2}$.  Then by the proof of Theorem \ref{s5t1}, it is not hard
to show that in this case we can still have the results of sections
\S3, \S4 and \S5 for the restriction of $\tilde{h}(t)$ to any $N_j$,
thus proving Theorem \ref{s5t2}.  We now proceed to show the above
injectivity radius bound.

 We claim that each $N_j$ is noncompact. In fact, by the curvature
 assumption on $M$, there exists $u$ such that $\sqrt{-1}\p\ol\p
 u=Ric_M$; see \cite{NST}. Let $\wt u$ be the lift of $u$ to $\wt
 M$. Then $\sqrt{-1}\p\ol\p\wt u=Ric_{\wt M}$. In 
particular, $\wt u$ is strictly plurisubharmonic on each $N_j$. Hence $N_j$ is noncompact.

By the proof in \cite{CZ3}, p. 25-26,  one can conclude that for any
$t_0>0$, there is a $\delta>0$ such that $h(t)$ 
has positive sectional curvature for $t_0<t\le  t_0+\delta$ when
restricted to $N_j$. Using the result of Gromoll-Meyer 
as before and using the fact that the curvature of $N_j$ is bounded
above by $C_1t^{-1}$ by Theorem \ref{shitamthm}, one
 can conclude that the injectivity radius of  $\tilde{h}(t)$ on $N_j$ is
 bounded below by $C_1t^{1/2}$ for some constant $C_1>0$ 
independent of $t$, $t_0$ and $j$. From this we can conclude that the
injectivity radius of $h(t_0)$ on $N_j$ is bounded 
below by $C_1t^{1/2}$. Hence the injectivity radius of $h(t)$ on $\wt
M$ is bounded below by $Ct^{\frac12}$ for some
 constant $C>0$ independent of $t$.  This completes the proof of the Theorem.

\end{proof}

 \bibliographystyle{amsplain}

\end{document}